\RequirePackage{fix-cm}
\documentclass[smallextended]{svjour3}       
\smartqed  
\usepackage{graphicx}
\usepackage{amsmath}
\usepackage{ dsfont }
\usepackage{algorithm}
\usepackage{algorithmic}
\usepackage{xcolor,colortbl}
\usepackage{subfigure}
\usepackage{booktabs}
\usepackage{tabularx}
\usepackage{placeins}
\usepackage{multirow}
\usepackage{algorithmic}
\usepackage{algorithm}
\allowdisplaybreaks
\begin{document}

\title{A Method of Sequential Log-Convex Programming for Engineering Design
}

\titlerunning{SLCP for Engineering Design}        

\author{Cody Karcher \and Robert Haimes
}


\institute{C. Karcher \and R. Haimes \at
              Department of Aeronautics and Astronautics, Massachusetts Institute of Technology \\ 
              77 Massachusetts Avenue, Cambridge, MA 02139, USA \\
              \email{ckarcher@mit.edu, haimes@mit.edu} 
}

\date{Received: date / Accepted: date}

\maketitle
\begin{abstract}

A method of Sequential Log-Convex Programming (SLCP) is constructed that exploits the log-convex structure present in many engineering design problems.  The mathematical structure of Geometric Programming (GP) is combined with the ability of Sequential Quadratic Program (SQP) to accommodate a wide range of objective and constraint functions, resulting in a practical algorithm that can be adopted with little to no modification of existing design practices.  Three test problems are considered to demonstrate the SLCP algorithm, comparing it with SQP and the modified Logspace Sequential Quadratic Programming (LSQP).  In these cases, SLCP shows up to a 77\% reduction in number of iterations compared to SQP, and an 11\% reduction compared to LSQP.  The airfoil analysis code XFOIL is integrated into one of the case studies to show how SLCP can be used to evolve the fidelity of design problems that have initially been modeled as GP compatible.  Finally, a methodology for design based on GP and SLCP is briefly discussed.

%
\end{abstract}

\section{Optimization for Engineering Design}

Two of the defining tasks of the engineering profession are \textit{analysis} and \textit{design}.  Analysis is a process by which engineers obtain data to quantify the behavior of a part, component, or system.  A full scale build and test is always the preferred approach to obtaining the highest quality analysis data, but external factors like cost and schedule make it impossible to utilize these methods every time an analysis must be performed.  

As a result, analysis is often performed via \textit{simulation}, where a set of simplified physics is modeled mathematically and then used to obtain the necessary analysis data.  Given the abundance of computing resources in the modern era, simulation has become the standard method for performing engineering analysis.  Thousands of analysis models (often distributed as software or code packages) have been released, including models for computational fluid dynamics (CFD), finite element analysis (FEA), electromagnetic simulation, and thermal analysis, just to name a few.

In general these analysis models are highly complex, requiring great expertise both to develop and to run as a user.  Due to this complexity, analysis is often conceptualized and implemented as a \textit{black box}, where the user provides inputs and obtains outputs with no visibility into the actual simulation being run \cite{martins2013multidisciplinary}.  In this way, analysis tools can be thought of as functions:

\begin{equation}
    \mathbf{y} = f (\mathbf{x})
    \label{analysisBlackBox}
\end{equation}

\noindent
where $\mathbf{x}$ represents the inputs to the analysis and $\mathbf{y}$ represents the obtained outputs.  In essence, the vector $\mathbf{x}$ is the mathematical abstraction of the design of a particular engineering system, and vector $\mathbf{y}$ is how that system performs.

If the goal of analysis is to determine the performance $\mathbf{y}$ of some given system design $\mathbf{x}$, the \textit{design} process seeks to determine some vector $\mathbf{x}$ that satisfies performance criteria $\mathbf{y}$.  The first step in design is often to simplify the vector $\mathbf{x}$ to include only the most critical parameters that define the system.  Once complete, the vector $\mathbf{x}$ exists in a vector space $\mathds{R}^N$ (called the \textit{design space}) where $N$ is the number of design decisions that have been retained in $\mathbf{x}$, typically referred to as the \textit{design variables}.  Determining the appropriate values for these design variables is a challenging task because while analysis can be performed with a single ``function call'' to the black box method in Equation \ref{analysisBlackBox}, determining an appropriate design requires that many candidate designs be evaluated in order to determine the best one, $\mathbf{x}^*$.

Thus, design does not take the simple functional form of Equation \ref{analysisBlackBox}.  Instead, engineering design problems cast in the language of mathematics are optimization problems:
\begin{equation}
    \begin{aligned}
        & \underset{\mathbf{x}}{\text{minimize}}
        & & f(\mathbf{x})
    \end{aligned}
    \label{designOpt}
\end{equation}

\noindent
But engineering design problems rarely take the form of Equation \ref{designOpt}.  Instead, \textit{constraints} are often imposed on the problem:
\begin{equation}
    \begin{aligned}
        & \underset{\mathbf{x}}{\text{minimize}}
        & & f(\mathbf{x})\\
        & \text{subject to}
        & & \mathbf{g}_i (\mathbf{x}) \leq 0 , \; i = 1, \ldots, N\\
        &&& \mathbf{h}_j (\mathbf{x}) = 0 , \; j = 1, \ldots, M
    \end{aligned}
    \label{NLP}
\end{equation}

Constraints typically serve one of two purposes.  First is to impose an artificial limit on the system, such as a minimum dimension or a maximum cost.  Second is to represent a limit imposed by a fundamental law of physics, such as the maximum stress that can be carried in a material or the dynamics of Newton's Second Law.  

These physics based constraints are another way in which analysis models are often included in engineering design problems, and so functions $f(\mathbf{x})$, $\mathbf{g}_i (\mathbf{x})$, and $\mathbf{h}_j (\mathbf{x})$ must all be assumed in the general case to be highly complicated black box functions that are expensive to evaluate.  The optimization methods used most commonly in engineering design are therefore tailored to minimize the number of times the objective and constraint functions must be evaluated.  Of these methods, Sequential Quadratic Programming (SQP), Geometric Programming (GP), and Logspace Sequential Quadratic Programming (LSQP) are most relevant to this work.

\section{Foundations in Existing Optimization Algorithms}
\subsection{Sequential Quadratic Programming}
In general, the problem posed in Equation \ref{NLP} is referred to as a Non-Linear Program (NLP) and is difficult to solve.  Many algorithms exist for solving NLPs, but the Sequential Quadratic Programming (SQP) algorithm is highly effective and has been utilized across a wide range of scientific and engineering fields.  The SQP algorithm begins with an initial guess $\mathbf{x}_k$ and formulates a Quadratic Programming (QP) approximation of the NLP that is valid in the local region near $\mathbf{x}_k$.  The QP sub-problem takes the form \cite{boggs1996sequential,nocedal2006numerical,kraft1988software}:
\begin{equation}
    \begin{aligned}
        & \underset{\mathbf{d}}{\text{minimize}}
        & & f(\mathbf{x}_k) + 
            \nabla f(\mathbf{x}_k)^T \mathbf{d} + 
            \frac{1}{2} \mathbf{d}^T    \nabla^2 \mathcal{L} (\mathbf{x}_k)   \mathbf{d}\\
        & \text{subject to}
        & & g_i(\mathbf{x}_k) + \nabla g_i(\mathbf{x}_k)^T \mathbf{d} \leq 0 , \; i = 1, \ldots, N\\
        &&& h_j(\mathbf{x}_k) + \nabla h_j(\mathbf{x}_k)^T \mathbf{d}    = 0 , \; j = 1, \ldots, M \\
        &&& \mathbf{d} = \mathbf{x}-\mathbf{x}_k 
    \end{aligned}
    \label{SQP_subproblem}
\end{equation}

\noindent
so named because of the quadratic objective function.  Since QPs are known to be \textit{convex} \cite{boyd2009convex}, the optimization problem in Equation \ref{SQP_subproblem} can be used to reliably and efficiently produce a new guess $\mathbf{x}_{k+1}$.  The process can then be iterated until some convergence criteria is reached, returning the optimal solution to the original NLP, $\mathbf{x}^*$.

Together with interior point methods, SQP represents the current state of the art for solving constrained continuous non-linear optimization problems \cite{martins2021engineering} despite being nearly 60 years old \cite{boggs1996sequential}.  The widespread success of SQP can be traced back to two key attributes.  First, the QP sub-problem is easy to construct since the sub-problem only requires the function evaluations $f(\mathbf{x}_k)$, $g_i(\mathbf{x}_k)$, and $h_j(\mathbf{x}_k)$ and the gradients $\nabla f(\mathbf{x}_k)$, $\nabla g_i(\mathbf{x}_k)$, and $\nabla h_j(\mathbf{x}_k)$.  These quantities are generally simple to obtain regardless of the complexity of the true functions, making SQP applicable to a incredibly large number of problem formulations.  Second, the QP sub-problem is easy to solve because it is convex, established in great detail by Boyd and Vandenberghe \cite{boyd2009convex}.  This convexity is key as convex optimization problems can be solved reliably and efficiently, unlike most other NLPs.

In fact, the QP sub-problem is the \textit{best possible} convex approximation of the true NLP that can be derived from a Taylor series decomposition of functions $f(\mathbf{x})$, $\mathbf{g}_i (\mathbf{x})$, and $\mathbf{h}_j (\mathbf{x})$, since taking any more terms in either the objective or constraint approximations would result in a non-convex sub-problem.  The desire for an accurate sub-problem should be relatively intuitive, since the number of iterations required to solve the NLP decreases as sub-problem accuracy increases\footnote{Consider as a thought experiment the extreme case where the sub-problem exactly represents the original NLP.  The solution would be obtained in only one iteration.}.  But many other forms of convex optimization problems exist, including Linear Programs (LP), Semi-Definite Programs (SDP), Second Order Cone Programs (SOCP), and some Quadratically Constrained Quadratic Programs (QCQP), among others.  If the role of the sub-problem is only to efficiently produce a new guess $x_{k+1}$, any one of these convex forms could easily be used in place of the QP sub-problem.  These theoretical approaches are generalized under the classification of Sequential Convex Programming (SCP) \cite{boyd2015sequential,duchi2018sequential}.

Of these SCP variations, only Sequential Quadratically Constrained Quadratic Programming (SQCQP) has received much attention in the literature \cite{anitescu2002superlinearly,tang2008sequential,liu2020method,jian2021qcqp}, but these methods struggle with complications in computing constraint curvature and in handling non-convex QCQP sub-problems\footnote{It is theoretically possible to model nearly any constrained continuous optimization as a QCQP regardless of whether convex structure exists or not, significantly limiting the general usefulness of the form.}.  Why other forms of SCP have not been studied is not clear, but any method of SCP should abide by the following criteria:
\begin{enumerate}
    \item Ease of construction should be comparable to the QP sub-problem of SQP (ie, use only $f(\mathbf{x}_k)$, $g_i(\mathbf{x}_k)$, and $h_j(\mathbf{x}_k)$ and the gradients $\nabla f(\mathbf{x}_k)$, $\nabla g_i(\mathbf{x}_k)$, and $\nabla h_j(\mathbf{x}_k)$ in sub-problem construction) 
    \item Exhibits a convex structure, and therefore easily solved
    \item Captures the underlying NLP more accurately than the QP sub-problem of SQP
\end{enumerate}

An LP based algorithm would rarely be superior to SQP, algorithms based on SDP or SOCP do not have an obvious construction method for the sub-problem, and the lack of convexity in some QCQPs has already been discussed.  So, is it possible to develop a method of SCP that satisfies all three criteria?  To answer this question one key building block remains, as recent literature has suggested a clear front runner for the type of convex optimization that should be used for engineering design: Geometric Programming.
\subsection{Geometric Programming}
A Geometric Program (GP) is a specific type of optimization formulation built from two classes of functions: monomials and posynomials.  A monomial function is defined as the product of a leading constant with each variable raised to a real power \cite{boyd2007tutorial}:
\begin{equation}
    m(\textbf{x}) = c{x_1}^{a_1}{x_2}^{a_2}...\;{x_n}^{a_n} = c \prod_{i=1}^{N} x_i^{a_i}
\end{equation}
A posynomial is simply the sum of monomials \cite{boyd2007tutorial}, which can be defined in notation as:
\begin{equation}
  p(\textbf{x}) = m_1(\textbf{x}) +  m_2(\textbf{x}) + ... + m_n(\textbf{x}) = \sum_{k=1}^{K} c_k \prod_{i=1}^{N} x_i^{a_{ik}}
\end{equation}
From these two building blocks, it is possible to construct the definition of a GP in standard form \cite{boyd2007tutorial}:
\begin{equation}
    \begin{aligned}
        & \underset{\mathbf{x}}{\text{minimize}}
        & & p_0 (\textbf{x}) \\
        & \text{subject to}
        & & m_i (\textbf{x}) = 1, \; i = 1, \ldots, N\\
        &&& p_j (\textbf{x}) \leq 1, \; j = 1, \ldots, M
    \end{aligned}
    \label{GP_standard}
\end{equation}

\noindent
If the general NLP (Equation \ref{NLP}) can be written in GP standard form (Equation \ref{GP_standard}) then it can be solved with great efficiency, since upon log transformation\footnote{In some of the literature, the transformation considered in this work is referred to as a log-log transformation since both dependent and independent variables are transformed. In all cases here, logspace, log-convexity, log transformation etcetera could equivalently be called log-log space, log-log convexity, log-log transformation and similar.} geometric programs become \textit{convex} \cite{boyd2007tutorial}.

The advantage of the GP form is that it is far more representative of many engineering design problems than the QP formulation.  The benefits of GP for engineering design has been well established in the literature \cite{clasen1984solution,greenberg1995mathematical,boyd2005digital,boyd2001optimal,li2004robust,xu2004oracle,jabr2005application,chiang2005geometric,chiang2007power,kandukuri2002optimal,marin2007optimization,vera2010optimization,preciado2014optimal,misra2014optimal,sela2015control} (see the original compilation in \cite{agrawal2019disciplined}), and has seen specific benefit for aircraft design \cite{hoburg2014geometric,torenbeek2013advanced,hoburg2013fast,kirschen2016signomial,brown2018vehicle,york2018efficient,burton2018solar,lin2020simultaneous,kirschen2018application,york2018turbofan,saab2018robust,hall2018assessment}.  Given that geometric programs are often more accurate at modeling engineering design problems, and that they can be solved just as efficiently as the quadratic programs that form the core of the SQP algorithm, GPs are a strong candidate for use in Sequential Convex Programming.

\subsection{Logspace Sequential Quadratic Programming}
The first attempts to leverage geometric programming in a sequential optimization algorithm were simply applications of SQP under the log transformation that makes GPs convex \cite{kirschen2018power,karcher2021logspace}.  Consider a slight modification of the general NLP \cite{karcher2021logspace}:
\begin{equation}
    \begin{aligned}
        & \underset{\mathbf{x}}{\text{minimize}}
        & & f(\mathbf{x})\\
        & \text{subject to}
        & & \mathbf{g}_i (\mathbf{x}) \leq 1 , \; i = 1, \ldots, N\\
        &&& \mathbf{h}_j (\mathbf{x}) = 1 , \; j = 1, \ldots, M
    \end{aligned}
    \label{nlpsf}
\end{equation}

Under the GP transformation $y_i = \log{x_i}$, or equivalently $x_i = e^{y_i}$, the problem becomes:
\begin{equation}
    \begin{aligned}
        & \underset{\mathbf{y}}{\text{minimize}}
        & & \log{ f (e^{\textbf{y}}) } \\
        & \text{subject to}
        & & \log{ \mathbf{g}_i (e^{\textbf{y}}) } \leq 0, \; i = 1, \ldots, N\\
        &&& \log{ \mathbf{h}_j (e^{\textbf{y}}) } = 0, \; j = 1, \ldots, M
    \end{aligned}
    \label{NLP_logtransformed}
\end{equation}

\noindent
which makes the new QP sub-problem \cite{karcher2021logspace}:
\begin{equation}
    \begin{aligned}
        & \underset{\mathbf{d}}{\text{minimize}}
        & & \log{f(\mathbf{x}_k)} + 
            \frac{1}{f(\mathbf{x}_k)}\left( \mathbf{x}_k \odot \nabla f(\mathbf{x}_k) \right)^T \mathbf{d} + 
            \frac{1}{2} \mathbf{d}^T  \nabla^2 \mathcal{L}(\mathbf{y}_k)   \mathbf{d}\\
        & \text{subject to}
        & & \log{g_i(\mathbf{x}_k)} + \frac{1}{g_i(\mathbf{x}_k)}\left( \mathbf{x}_k \odot \nabla g_i(\mathbf{x}_k) \right)^T \mathbf{d} \leq 0, \; i = 1, \ldots, N\\
        &&& \log{h_j(\mathbf{x}_k)} + \frac{1}{h_j(\mathbf{x}_k)}\left( \mathbf{x}_k \odot \nabla h_j(\mathbf{x}_k) \right)^T \mathbf{d}    = 0 , \; j = 1, \ldots, M \\
        &&& \mathbf{d} = \mathbf{y}-\log{\mathbf{x}_k} \\
        &&& \mathbf{y} = \log{\mathbf{x}}
    \end{aligned}
    \label{lsqp_sp}
\end{equation}

This approach showed great promise when applied to a series of engineering design case studies \cite{karcher2021logspace}, but can be taken one step further.  The sub-problem defined by Equation \ref{lsqp_sp} represents monomial constraints exactly, but gives no consideration to posynomial functions.  Under transformation, a posynomial constraint becomes:
\begin{equation}
    \log \left( \sum_j \exp \left( P_j (\mathbf{d} + \log{\mathbf{x}_k})  + q_j \right) \right) \leq 0
    \label{posyLog}
\end{equation}
which if inserted into Equation \ref{lsqp_sp} \textit{leaves the sub-problem convex}.  Since convex problems can be readily solved \cite{boyd2009convex} it is possible to model these posynomial constraints directly, without resorting to the linearized form.  

\section{The Mathematics of SLCP}
\subsection{Mathematical Definition of the SLCP Method}
Consider that rather than representing the general non-linear program with Equation \ref{nlpsf}, the constraints that are GP compatible (posynomials and monomials) are given special treatment:
\begin{equation}
    \begin{aligned}
        & \underset{\mathbf{x}}{\text{minimize}}
        & & f(\mathbf{x})\\
        & \text{subject to}
        & & \mathbf{p} (\mathbf{x}) \leq 1 \\
        &&& \mathbf{m} (\mathbf{x}) = 1 \\
        &&& \mathbf{g} (\mathbf{x}) \leq 1 \\
        &&& \mathbf{h} (\mathbf{x}) = 1 
    \end{aligned}
    \label{slcpStructure}
\end{equation}

The constraints $\mathbf{g} (\mathbf{x}) \leq 1$ and $\mathbf{h} (\mathbf{x}) = 1$ will be linearized in the sub-problem just as in LSQP \cite{karcher2021logspace}, but posynomials and monomials can be transformed and imposed directly in the sub-problem.  Following this procedure yields the following sub-problem form:
\begin{equation}
    \begin{aligned}
        & \underset{\mathbf{d}}{\text{minimize}}
       & & \log{f(\mathbf{x}_k)} + 
            \frac{1}{f(\mathbf{x}_k)}\left( \mathbf{x}_k \odot \nabla f(\mathbf{x}_k) \right)^T \mathbf{d} + 
            \frac{1}{2} \mathbf{d}^T  \nabla^2 \mathcal{L}_R (\mathbf{y}_k)   \mathbf{d}\\
        & \text{subject to}
        & & \log \left( \sum_j \exp \left( P_j (\mathbf{d} + \log{\mathbf{x}_k})  + q_j \right) \right) \leq 0\\
        &&& A_m (\mathbf{d} + \log{\mathbf{x}_k}) + b_m \leq 0 \\
        &&& \log{g(\mathbf{x}_k)} + \frac{1}{g(\mathbf{x}_k)}\left( \mathbf{x}_k \odot \nabla g(\mathbf{x}_k) \right)^T \mathbf{d} \leq 0 \\
        &&& \log{h(\mathbf{x}_k)} + \frac{1}{h(\mathbf{x}_k)}\left( \mathbf{x}_k \odot \nabla h(\mathbf{x}_k) \right)^T \mathbf{d}    = 0 \\
        &&& \mathbf{d} = \mathbf{y}-\log{\mathbf{x}_k} \\
        &&& \mathbf{y} = \log{\mathbf{x}}
    \end{aligned}
    \label{slcpSubproblem_nbb}
\end{equation}

Solving the general non-linear program in Equation \ref{slcpStructure} via a series of sub-problems defined by Equation \ref{slcpSubproblem_nbb} is proposed here as a method of of Sequential Log-Convex Programming (SLCP), since the quadratic programing sub-problem of SQP has now been replaced with a log-convex programming (LCP) sub-problem.  
\subsection{The Reduced Lagrangian}
Though the primary distinction between the SLCP method proposed here and the LSQP algorithm in the literature \cite{karcher2021logspace} is in the handling of posynomial constraints, the objective function also requires a minor update.  The LSQP sub-problem inherits its quadratic objective function directly from SQP, which utilizes the Hessian of the Lagrangian function, defined as:
\begin{equation}
    \mathcal{L}(\mathbf{y},\lambda) = \log{f(\mathbf{x})} + \lambda \log{\mathbf{p}(\mathbf{x}_k)} 
                                                          + \lambda \log{\mathbf{m}(\mathbf{x}_k)} 
                                                          + \lambda \log{\mathbf{g}(\mathbf{x}_k)} 
                                                          + \lambda \log{\mathbf{h}(\mathbf{x}_k)} 
\end{equation}

The primary purpose of using the $\nabla^2 \mathcal{L} (\mathbf{x}_k)$ rather than $\nabla^2 f (\mathbf{x}_k)$ is to include some second order information from the constraints in the sub-problem \cite{boggs1996sequential,nocedal2006numerical}.  In the case of SLCP, some of the constraints with higher order curvature are now being represented directly, and so attempting to approximate the second order information of these constraints in the objective function causes a conflict between the approximated curvature and the true curvature that is now being fully captured.  Thus, those constraints must be left out of the second order Hessian approximation.  

The SLCP algorithm therefore utilizes a Reduced Lagrangian, which does not include the constraints which are exactly represented:
\begin{equation}
    \mathcal{L}_R(\mathbf{y},\lambda) = \log{f(\mathbf{x})} + \lambda \log{\mathbf{g}(\mathbf{x}_k)} 
                                                            + \lambda \log{\mathbf{h}(\mathbf{x}_k)} 
    \label{reducedLagrangian}
\end{equation}

The use of this Reduced Lagrangian is critical to the success of the algorithm.  Imposing exact constraints without this modification performs worse than strict LSQP.

\subsection{Limitations and Potential Improvements}
This proposed method of SLCP shares many of the same drawbacks as LSQP \cite{karcher2021logspace}.  Due to the log transformation, variables must be strictly positive, and remain strictly positive during the solve.  Since it is not possible to have a negative mass, length, volume, or similar, this limitation proves remarkably non-intrusive for many engineering design problems.

Likewise, functions $f(\mathbf{x})$, $\mathbf{g}_i (\mathbf{x})$, and $\mathbf{h}_j (\mathbf{x})$ must be positive at the initial guess and stay positive throughout the solution.  This concern is addressed in more depth in previous work \cite{karcher2021logspace}, but essentially, it is possible to mitigate this concern through intelligent function construction, and the step size can always be constrained to ensure these functions remain positive.  

One possible area for future improvement is the use of the quadratic objective function.  Some effort here was given to replacing the quadratic objective with a GP compatible posynomial function, but utilizing the quadratic objective with the BFGS approximation provided the best result and so was carried over from LSQP (with the Reduced Lagrangian modification).  A deep dive into modifying BFGS to accommodate a non-quadratic objective was viewed as being beyond scope.

Another compelling reason for keeping the quadratic objective in this work was that the sub-problem proposed in Equation \ref{slcpSubproblem_nbb} reverts to the LSQP sub-problem in the absence of posynomial constraints, which provides a cornerstone for comparison and debugging.  And since LSQP is simply an application of SQP in log transformed space \cite{karcher2021logspace},  the body of literature surrounding SQP could still be utilized with only minimal modification.

A second potential improvement stems from the realization that the SLCP sub-problem defined by Equation \ref{slcpSubproblem_nbb} is not the only possible SLCP sub-problem.  Indeed, work by Agrawal et. al. \cite{agrawal2019disciplined} on Disciplined Geometric Programming suggests that there are families of functions beyond posynomials that are log-convex and could therefore receive similar treatment in Equation \ref{slcpSubproblem_nbb}.  Developing a more general SLCP method that accounts for these functions will be the subject of future work, but was determined to be out of scope for this paper, especially given the success of GP form in engineering design without these additional functions.

\section{An Algorithm for the Proposed SLCP Method} \label{slcpAlgo}
Algorithm \ref{alg:slcp} outlines a method for implementing the proposed SLCP method programmatically.  Though Algorithm \ref{alg:slcp} is very similar to the well known SQP algorithm and the LSQP algorithm outlined in the literature \cite{karcher2021logspace}, a few minor changes should be noted.  

Similar to LSQP, Algorithm \ref{alg:slcp} utilizes the gradients $\frac{\partial \log{ f (e^{\textbf{y}}) }}{\partial y_i}$, which can be directly computed from the original gradients $\frac{\partial f}{\partial x_i}$ for minimal computational expense \cite{boyd2007tutorial,karcher2021logspace}:
\begin{equation}
    \frac{\partial \log{ f (e^{\textbf{y}}) }}{\partial y_i} = \frac{x_i}{f(\mathbf{x})} \frac{\partial f}{\partial x_i} 
    \label{gradient_F}
\end{equation}

The use of the Reduced Lagrangian also necessitates modification to the damped BFGS method \cite{nocedal2006numerical}:
\begin{equation}
\begin{aligned}
    \mathbf{B}_{k+1} & = \mathbf{B}_k - \frac{\mathbf{B}_k s_k s_k^T \mathbf{B}_k}{s_k^T \mathbf{B}_k s_k} + \frac{r_k r_k^T}{s_k^T r_k} \\
    s_k & = \alpha_k \mathbf{d}_x \\
    z_k & = \nabla \mathcal{L}_{R_{k+1}}(f,g,h,\mathbf{y}_{k+1},\mu_{k+1}) - \nabla \mathcal{L}_{R_{k}}(f,g,h,\mathbf{y}_{k},\mu_{k+1})  \\
    r_k & = \theta_k z_k + (1-\theta_k) \mathbf{B}_k s_k \\
    \theta_k & =\left\{
      \begin{array}{@{}ll@{}}
        1 & \text{if}\ s_k^T z_k \geq 0.2 s_k^T \mathbf{B}_k s_k \\
        (0.8s_k^T \mathbf{B}_k s_k)/(s_k^T \mathbf{B}_k s_k - s_k^T z_k ) & \text{if}\ s_k^T z_k <    0.2 s_k^T \mathbf{B}_k s_k   
      \end{array}\right.
\end{aligned}
\label{modifiedBFGS}
\end{equation}

Finally, the sub-problem is relaxed here to handle the problem of inconsistent constraints frequently faced by SQP \cite{nocedal2006numerical}:
\begin{equation}
    \begin{aligned}
       & \underset{\mathbf{d}}{\text{minimize}}
       & & \log{f(\mathbf{x}_k)} + 
            \frac{1}{f(\mathbf{x}_k)}\left( \mathbf{x}_k \odot \nabla f(\mathbf{x}_k) \right)^T \mathbf{d} + 
            \frac{1}{2} \mathbf{d}^T  \nabla^2 \mathcal{L}_R (\mathbf{y}_k)   \mathbf{d} + K \sum_{i=1}^{N_{con}} \sigma_i^2 \\
        & \text{subject to}
        & & \log \left( \sum_j \exp \left( P_j (\mathbf{d} + \log{\mathbf{x}_k})  + q_j \right) \right) \leq \sigma_i\\
        &&& A_m (\mathbf{d} + \log{\mathbf{x}_k}) + b_m \leq \sigma_i \\
        &&& \log{g(\mathbf{x}_k)} + \frac{1}{g(\mathbf{x}_k)}\left( \mathbf{x}_k \odot \nabla g(\mathbf{x}_k) \right)^T \mathbf{d} \leq \sigma_i \\
        &&& \log{h(\mathbf{x}_k)} + \frac{1}{h(\mathbf{x}_k)}\left( \mathbf{x}_k \odot \nabla h(\mathbf{x}_k) \right)^T \mathbf{d}    = \sigma_i \\
        &&& \mathbf{d} = \mathbf{y}-\log{\mathbf{x}_k} \\
        &&& \mathbf{y} = \log{\mathbf{x}}
    \end{aligned}
    \label{slcp_sp_relaxed}
\end{equation}

\begin{algorithm}
\caption{Sequential Log-Convex Programming (SLCP)}
\label{alg:slcp}
\begin{algorithmic}[1]
\STATE{Given $\mathbf{x}_0$}
\STATE{Construct standard form}
\STATE{Compute $\mathbf{y}_0 = \log ( \mathbf{x}_0 )$}
\STATE{Initialize logspace Lagrange multipliers, $\mu_0 \leftarrow \mathbf{1}$}
\STATE{Initialize the matrix $\mathbf{B} \leftarrow \mathbf{I}$ the approximation of $\nabla^2 \mathcal{L}_R(\mathbf{y},\mu)$ \cite{nocedal2006numerical}}
\STATE{Compute $f(\mathbf{x}_0)$, $\mathbf{g}(\mathbf{x}_0)$, $\mathbf{h}(\mathbf{x}_0)$, $\mathbf{p}(\mathbf{x}_0)$, and $\mathbf{m}(\mathbf{x}_0)$ }
\STATE{Compute $\log f(\mathbf{x}_0)$, $\log \mathbf{g}(\mathbf{x}_0)$, $\log \mathbf{h}(\mathbf{x}_0)$, $\log \mathbf{p}(\mathbf{x}_0)$, and $\log \mathbf{m}(\mathbf{x}_0)$}
\STATE{Compute $\nabla f(\mathbf{x}_0)$, $\nabla \mathbf{g}(\mathbf{x}_0)$, $\nabla \mathbf{h}(\mathbf{x}_0)$, $\nabla \mathbf{p}(\mathbf{x}_0)$, and $\nabla \mathbf{m}(\mathbf{x}_0)$}
\STATE{Compute $\nabla \log f(e^{\mathbf{y}_0})$, $\nabla \log \mathbf{g}(e^{\mathbf{y}_0})$, $\nabla \log \mathbf{h}(e^{\mathbf{y}_0})$, $\nabla \log \mathbf{p}(e^{\mathbf{y}_0})$, and $\nabla \log \mathbf{m}(e^{\mathbf{y}_0})$ via Equation \ref{gradient_F}}
\FOR{$k=0$ to maxIter}
  \STATE{Solve the convex sub-problem (Equation \ref{slcp_sp_relaxed}) to obtain $\mathbf{d}_y$ and $\mathbf{d}_{\mu}$ \cite{nocedal2006numerical}}
  \STATE{Compute the step size $\alpha_k$ via inexact line search \cite{nocedal2006numerical,matlab2020constrained}}
  \STATE{$\mathbf{y}_{k+1} \leftarrow \mathbf{y}_k + \alpha_k \mathbf{d}_y$}
  \STATE{$\mathbf{x}_{k+1} \leftarrow \exp(\mathbf{y}_{k+1})$}
  \STATE{$\mu_{k+1} \leftarrow \mu_k + \alpha_k \mathbf{d}_{\mu}$}
  \STATE{Compute $f(\mathbf{x}_{k+1})$, $\mathbf{g}(\mathbf{x}_{k+1})$, $\mathbf{h}(\mathbf{x}_{k+1})$, $\mathbf{p}(\mathbf{x}_{k+1})$, and $\mathbf{m}(\mathbf{x}_{k+1})$ }
  \STATE{Compute $\log f(\mathbf{x}_{k+1})$, $\log \mathbf{g}(\mathbf{x}_{k+1})$, $\log \mathbf{h}(\mathbf{x}_{k+1})$, $\log \mathbf{p}(\mathbf{x}_{k+1})$, and $\log \mathbf{m}(\mathbf{x}_{k+1})$}
  \STATE{Compute $\nabla f(\mathbf{x}_{k+1})$, $\nabla \mathbf{g}(\mathbf{x}_{k+1})$, $\nabla \mathbf{h}(\mathbf{x}_{k+1})$, $\nabla \mathbf{p}(\mathbf{x}_{k+1})$, and $\nabla \mathbf{m}(\mathbf{x}_{k+1})$}
  \STATE{Compute $\nabla \log f(e^{\mathbf{y}_{k+1}})$, $\nabla \log \mathbf{g}(e^{\mathbf{y}_{k+1}})$, $\nabla \log \mathbf{h}(e^{\mathbf{y}_{k+1}})$, $\nabla \log \mathbf{p}(e^{\mathbf{y}_{k+1}})$, and $\nabla \log \mathbf{m}(e^{\mathbf{y}_{k+1}})$ via Equation \ref{gradient_F}}
  \STATE{Compute $\nabla \mathcal{L}_k(\mathbf{y}_k,\mu_{k+1})$ \cite{nocedal2006numerical}}
  \STATE{Compute $\nabla \mathcal{L}_{k+1}(\mathbf{y}_{k+1},\mu_{k+1})$ \cite{nocedal2006numerical}}
  \STATE{Compute $\nabla \mathcal{L}_{R_k}(\mathbf{y}_k,\mu_{k+1})$ \cite{nocedal2006numerical}}
  \STATE{Compute $\nabla \mathcal{L}_{R_{k+1}}(\mathbf{y}_{k+1},\mu_{k+1})$ \cite{nocedal2006numerical}}
  \IF{$\nabla \mathcal{L}_{k+1}(\mathbf{y}_{k+1},\mu_{k+1}) < \varepsilon_{GL}$ \cite{nocedal2006numerical}}
    \RETURN{$\mathbf{x}_{k+1}$}
  \ELSIF{$|| \mathbf{d}_{x} || < \varepsilon_{dx}$ \cite{matlab2020constrained}}
    \RETURN{$\mathbf{x}_{k+1}$}
  \ELSE
    \STATE{Perform a modified damped BFGS update on matrix $\mathbf{B}$ (Equation \ref{modifiedBFGS}) \cite{nocedal2006numerical}}
    \STATE{$k \leftarrow k+1$}
  \ENDIF
\ENDFOR
\RETURN{$x_k$, maximum iteration count reached}
\end{algorithmic}
\end{algorithm}

One important note is that unlike LSQP, this SLCP algorithm cannot utilize existing SQP solvers due to the fundamentally different nature of the sub-problem construction.  In the work presented here, Algorithm \ref{alg:slcp} is implemented in a custom python suite, but the sub-problems are solved using the appropriate CVXOPT \cite{andersen2013cvxopt} solver.
\section{A Simple Test Case}
Before launching into complex design examples, it is valuable to gain intuition through the use of a simple example.  Consider the following geometric program:
\begin{equation}
    \begin{aligned}
        & \underset{\mathbf{x,y}}{\text{minimize}}
        & & \frac{1}{x^{0.1}} + 15 x^{0.01} + \frac{1}{y^{0.1}} + 15 y^{0.01} \\
        & \text{subject to}
        & & 0.01 x^{-1.1} + x^{0.1}  + y \leq 1 \\
        &&& x \geq \epsilon \\
        &&& y \geq \epsilon \\
    \end{aligned}
    \label{exampleProblem}
\end{equation}

\noindent
where $\epsilon$ is some small positive value, in this case $10^{-9}$.  Under the log transformation utilized by GP, LSQP, and SLCP, this problem becomes:

\begin{equation}
    \begin{aligned}
        & \underset{u,v}{\text{minimize}}
        & & \log \left ( e^{-0.1u} + 15e^{0.01u} + e^{-0.1v} + 15e^{0.01v}\right ) \\
        & \text{subject to}
        & & \log \left ( e^v + 0.01e^{-1.1u} + e^{0.1u} \right ) \leq 0 \\
    \end{aligned}
    \label{exampleProblemTransformed}
\end{equation}

\begin{figure}[htb]
\centering     
\subfigure[As written (Equation \ref{exampleProblem})]{\includegraphics[width=0.48\textwidth]{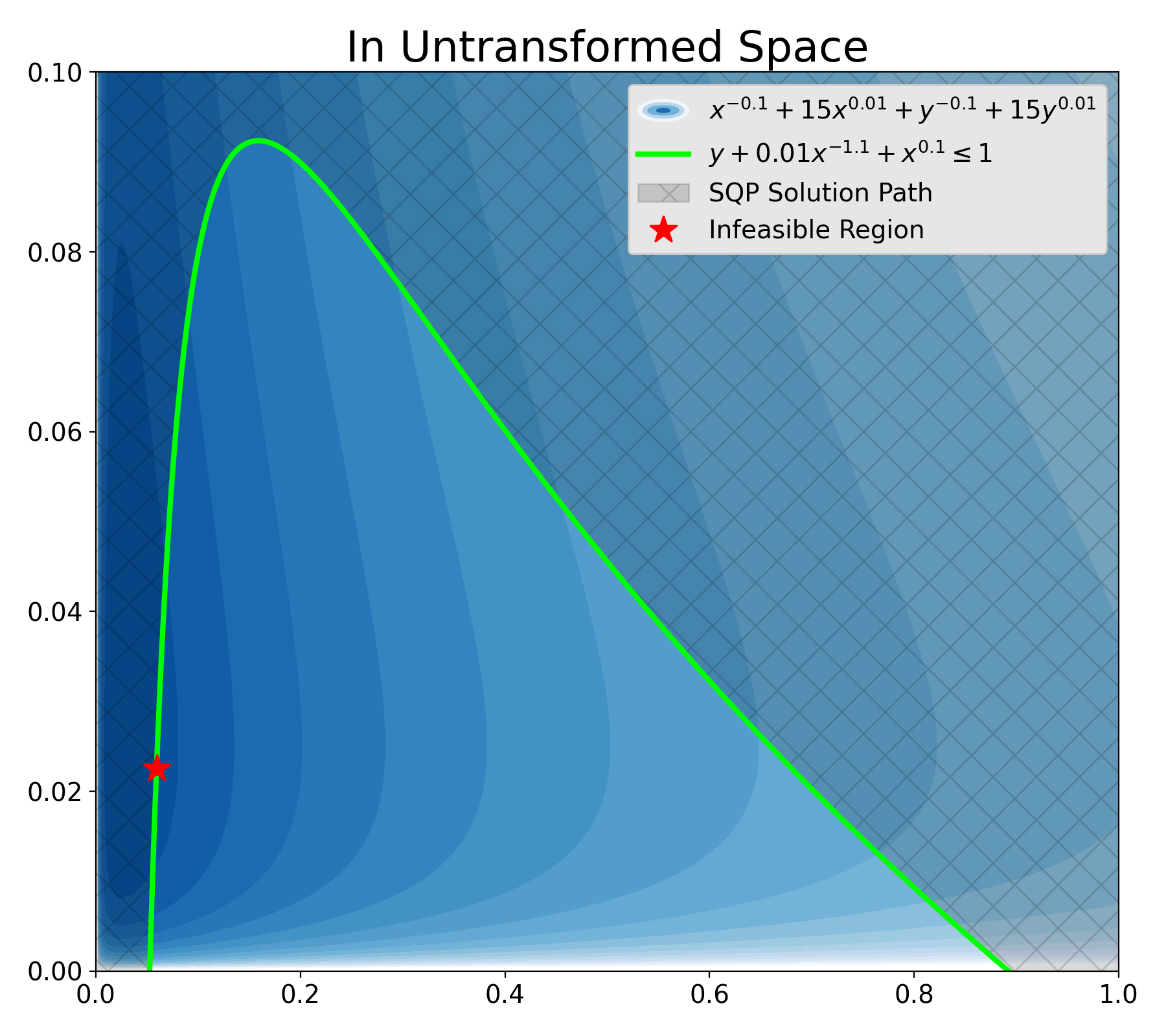} \label{seut}}
\subfigure[Under log transformation (Equation \ref{exampleProblemTransformed})]{\includegraphics[width=0.48\textwidth]{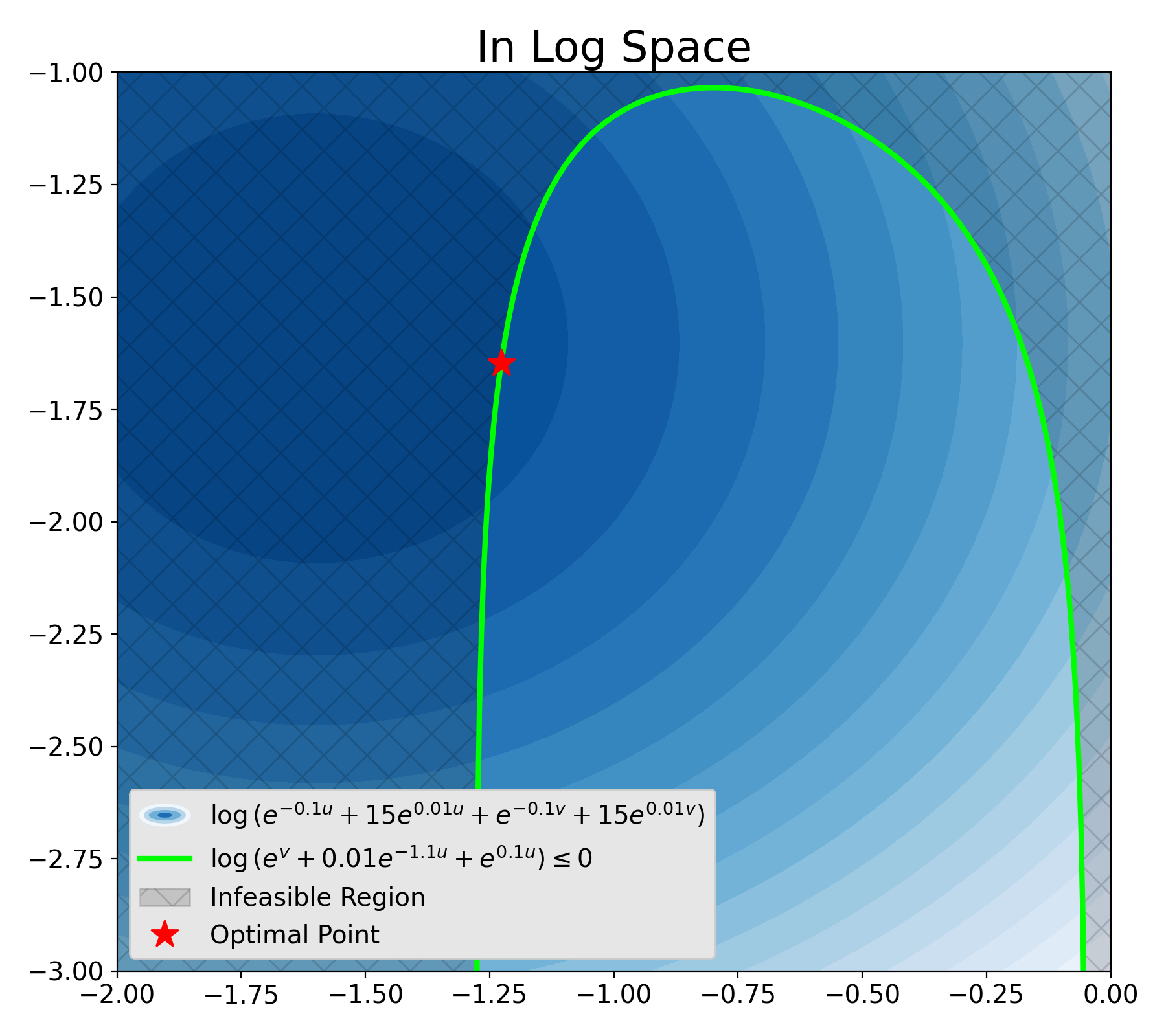} \label{setr}}
\caption{Visualizing the simple example problem in both the untransformed (Figure \ref{seut}) and log transformed (Figure \ref{setr}) spaces}
\label{exampleProblemPlot}
\end{figure}

The problem can be directly visualized, both in the original space and in the log transformed space, as seen in Figure \ref{exampleProblemPlot}.
\FloatBarrier

Figure \ref{exampleProblemPlot} highlights the clear advantage of the log transformation.  Figure \ref{seut} is characterized by long, skinny, irregular objective contours and a non-convex constraint, both of which make it poorly conditioned for solution with gradient based methods.  In contrast, Figure \ref{setr} has objective contours that are nearly circular in shape, and a constraint that carves out a convex set\footnote{The reader should not be confused by the use of the term `convex set' here because the constraint function is in fact a concave function.  A convex optimization problem is by definition the minimization of a convex objective function over a convex set generated by the constraint set, which is the case here.} with no cusps or drastic changes of curvature, making it highly conducive to gradient based optimization.  While it is true that not all optimization formulations will benefit from this transformation (see discussion of the Rosenbrock problem in Karcher \cite{karcher2021logspace}), the abundance of literature showing the applicability of geometric programming to engineering design problems, along with the success of the LSQP algorithm, indicates this transformation is a useful tool in many cases of interest.  

\begin{figure}[htb]
\centering     
\subfigure[With LSQP]{\includegraphics[width=0.48\textwidth]{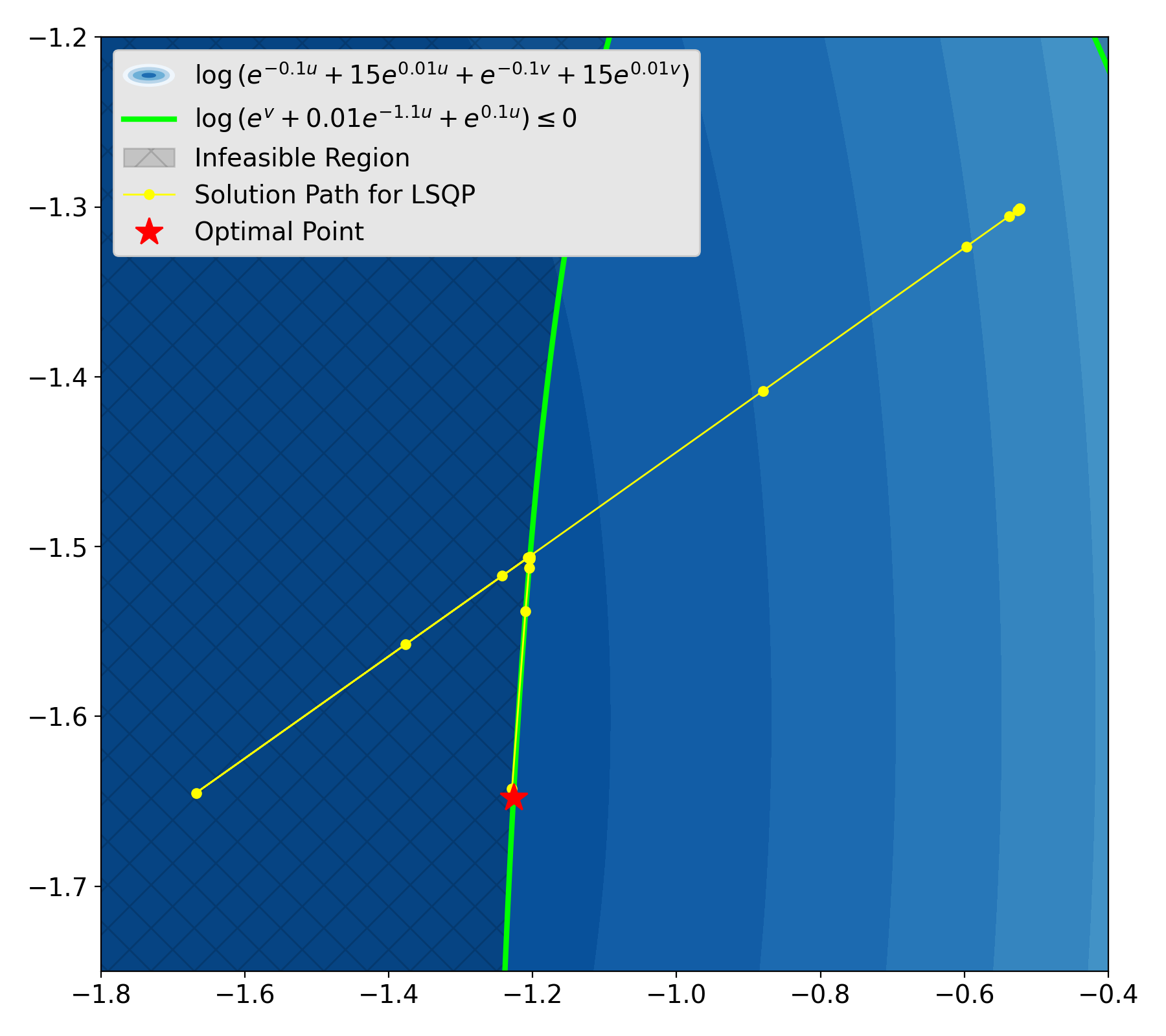} \label{lsqpPath}}
\subfigure[With SLCP]{\includegraphics[width=0.48\textwidth]{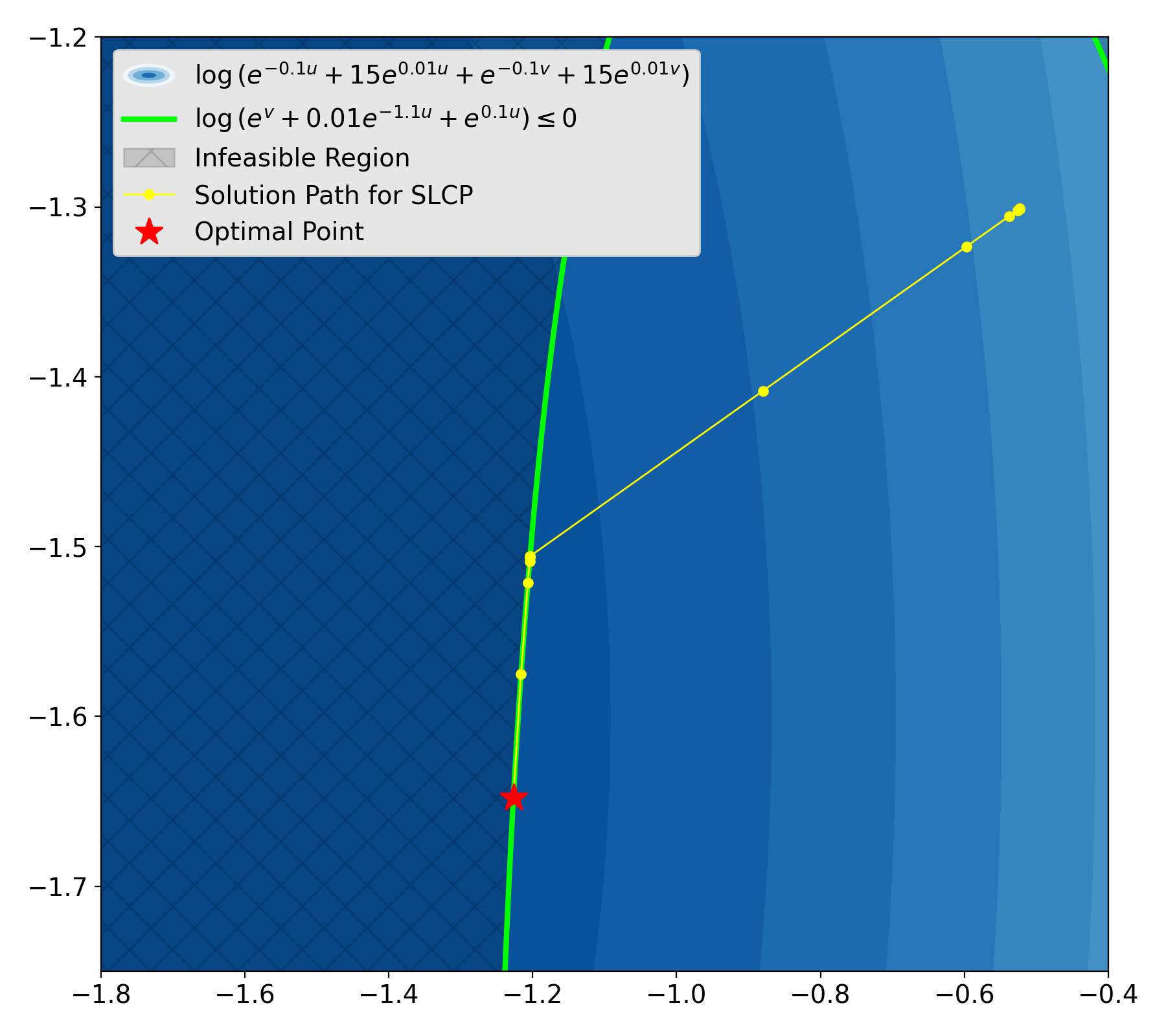} \label{slcpPath}}
\caption{Plotting the solution paths taken when solving the simple example problem with the LSQP (Figure \ref{lsqpPath}) and SLCP (Figure \ref{slcpPath}) algorithms}
\label{exampleProblemLSQPsolve}
\end{figure}

Now consider the difference between LSQP and SLCP.  Starting from $(x_0,y_0) = (0.3,0.05)$, Figure \ref{exampleProblemLSQPsolve} shows the two solution paths taken to reach the optimal solution (zoomed in to highlight the intermediate steps).  

\FloatBarrier

\begin{figure}[htb]
\centering     
\subfigure[The 6th iteration of LSQP]{\includegraphics[width=0.48\textwidth]{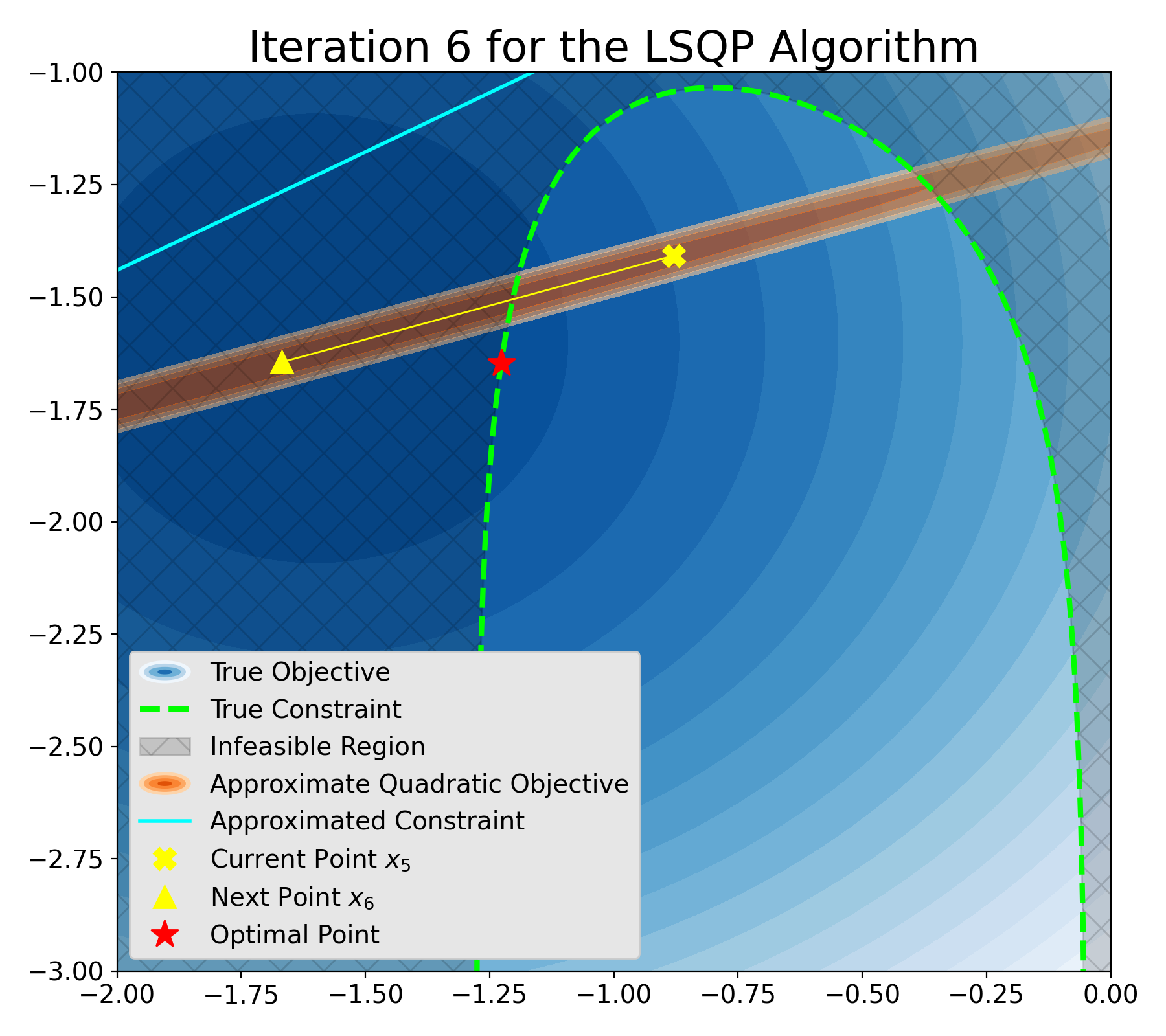} \label{itr6a}}
\subfigure[The 6th iteration of SLCP]{\includegraphics[width=0.48\textwidth]{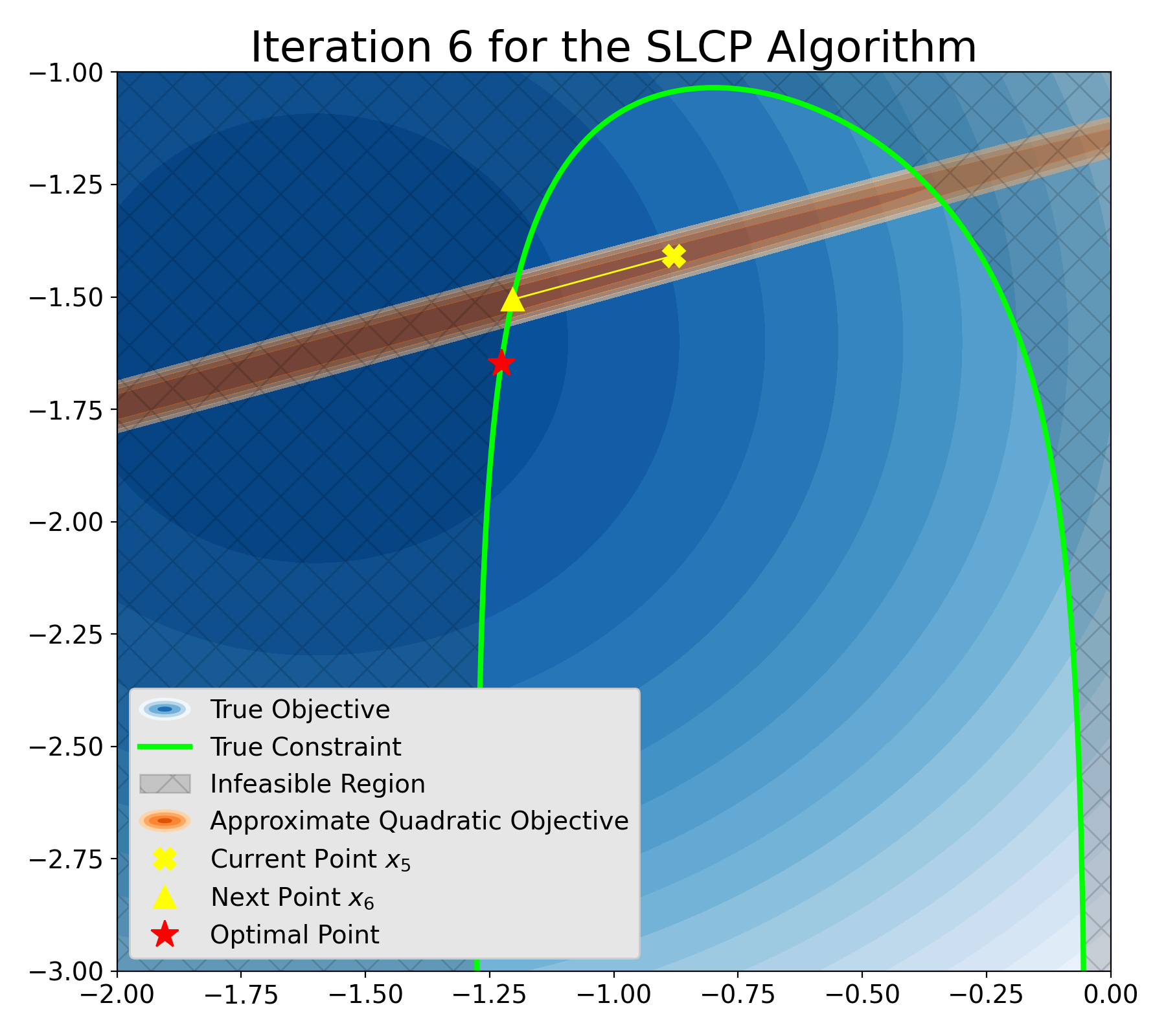} \label{itr6b}}
\caption{Comparing the critical 6th iteration of the LSQP and SLCP algorithms when solving the simple example problem}
\label{iteration6}
\end{figure}

Clearly, the LSQP algorithm overshoots the constraint at some point during the solution, and must work back into the feasible region before finally reaching convergence (Figure \ref{lsqpPath}).  It is the 6th iteration of both algorithms that distinguishes the performance, shown in Figure \ref{iteration6}.  

\FloatBarrier

The linear approximation in Figure \ref{itr6a} does not bound the QP sub-problem appropriately, resulting in an overshoot of the true constraint.  In contrast, Figure \ref{itr6b} shows the enforced posynomial bounding the step as desired, preventing overshoot.  In total, this enforcement saves the SLCP algorithm 4 iterations compared to LSQP.

\begin{figure}[htb]
\centering     
\subfigure[The final iteration of LSQP]{\includegraphics[width=0.48\textwidth]{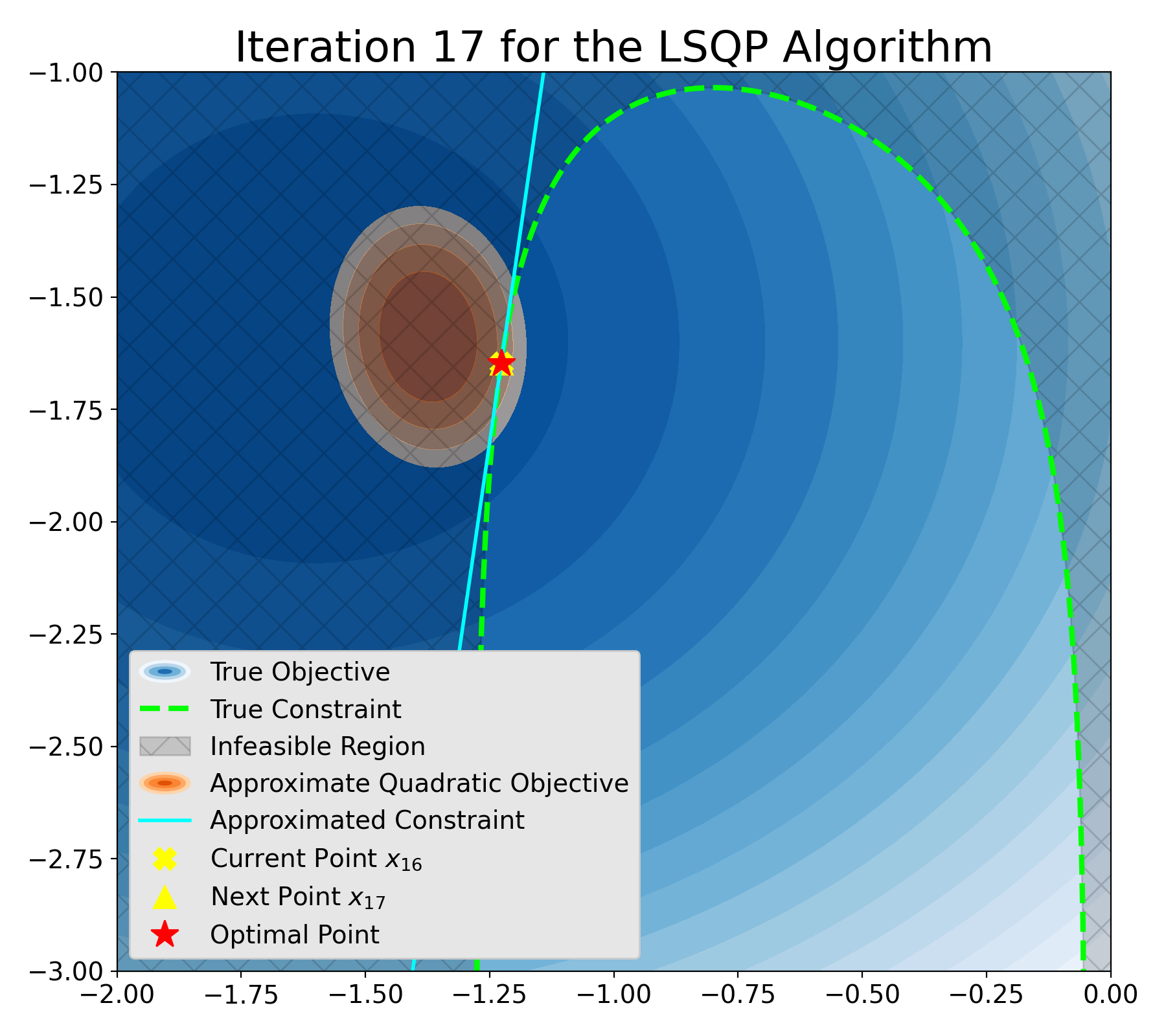} \label{itr17a}}
\subfigure[The final iteration of SLCP]{\includegraphics[width=0.48\textwidth]{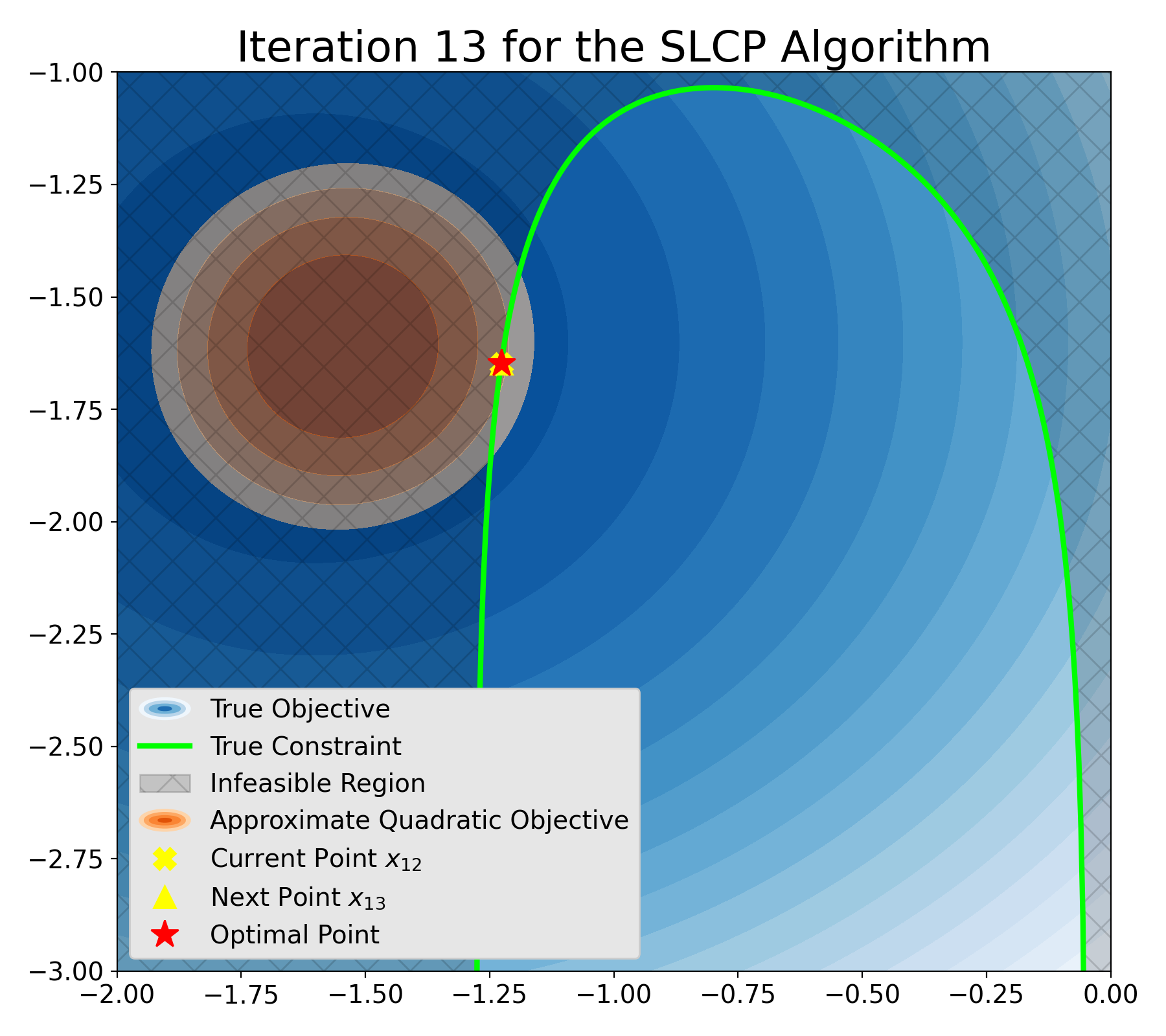} \label{itr13b}}
\caption{Comparing the second-order approximation of the objective function at the final iteration of the LSQP and SLCP algorithms}
\label{iterationFinal}
\end{figure}

The affect of the Reduced Lagrangian (Equation \ref{reducedLagrangian}) can also be visualized at the termination step of both algorithms, seen in Figure \ref{iterationFinal}.
\FloatBarrier

The approximated objective function in Figure \ref{itr13b} is a superior approximation of the true underlying objective function due to the absence of the constraint curvature terms that must be present in LSQP.

From this simple example, it is possible to draw the two conclusions that will be seen in the more extensive trials that follow.  First, SLCP will have its widest performance gap with LSQP when more posynomial constraints are present in the formulation.  Perhaps phrased more accurately, SLCP will outperform LSQP whenever it arises during the solution process that the linear approximation of a posynomial substantially deviates from the posynomial itself along the search direction.  Second, SLCP will outperform LSQP by a wider margin when the initial guess is farther from the true optimal solution, as it is expected that the deviation between the posynomials and their linear approximations will grow wider over larger distances.  These two trends will both be observed in the results below.
\section{Evaluating Algorithm Performance}
\subsection{Methodology}
To systematically test the effectiveness of the proposed SLCP method, three engineering design test problems were selected from the literature (Floudas \cite{floudas2013handbook}, Kirschen-Ozturk \cite{kirschen2018power}, and Hoburg \cite{hoburg2014geometric}) and solved using one of three methods:
\begin{enumerate} 
    \item A python implementation of SQP (Floudas and Kirschen-Ozturk only) 
    \item A python implementation of LSQP as described in \cite{karcher2021logspace}
    \item A python implementation of SLCP as described in Section \ref{slcpAlgo}
\end{enumerate}
\noindent
This methodology is similar to the one used to demonstrate the effectiveness of the LSQP algorithm \cite{karcher2021logspace}.  Note that three variations of the Hoburg problem were considered, as will be discussed below (see Sections \ref{hoburgProblem}, \ref{hoburgProblem1}, and \ref{hoburgProblem3}).

For each of the 12 problem/algorithm combinations, 3000 trials were run starting from a random initial starting point.  In 1000 of these cases, the initial guess was bounded to be within +/- 10\% of the known optimum, another 1000 were bound within +/- 50\% of the known optimum, and the final 1000 were bounded to be within +/-80\% of the known optimum.

For each set of 1000 trials, a curve was constructed showing the fraction of cases that had converged within a certain number of iterations (Figures \ref{floudas_results_fig}, \ref{ko_results_fig}, \ref{h0_results_fig}, \ref{h1_results_fig}, and \ref{h3_results_fig}).  In these plots the ideal algorithm would have a ``$\Gamma$-like'' shape, converging all cases in only one iteration, and so curves closest to the upper left of the graph represent the superior algorithms.

Due to the computational expense of 36000 trials, the computational resources of the MIT SuperCloud were utilized \cite{reuther2018interactive}, and a limit of 500 iterations was placed on all 3 algorithms.

\subsection{Floudas Problem}
Floudas \cite{floudas2013handbook} offers the following optimization for the design of a heat exchanger:
\begin{equation}
    \begin{aligned}
         \underset{x_1,...,x_8}{\text{minimize}}  \quad &  x_1 + x_2 + x_3 \\
         \text{subject to}  \quad & \frac{833.33252 x_4}{x_2 x_6} + \frac{100}{x_6} - \frac{83333.333}{x_1 x_6} \leq 1 \\
         & \frac{1250 x_5}{x_2 x_7} + \frac{x_4}{x_7} - \frac{1250 x_4}{x_2 x_7} \leq 1 \\
         & \frac{1250000}{x_3 x_8} + \frac{x_5}{x_8} - \frac{2500 x_5}{x_3 x_8}  \leq 1 \\
         & 0.0025 x_4 + 0.0025 x_6  \leq 1 \\
         & -0.0025 x_4 + 0.0025 x_5 + 0.0025 x_7  \leq 1 \\
         & -0.01 x_5 + 0.01 x_8  \leq 1
    \end{aligned}
    \label{floudas}
\end{equation}

\begin{figure}[htb]
\centering     
\includegraphics[width=0.9\textwidth]{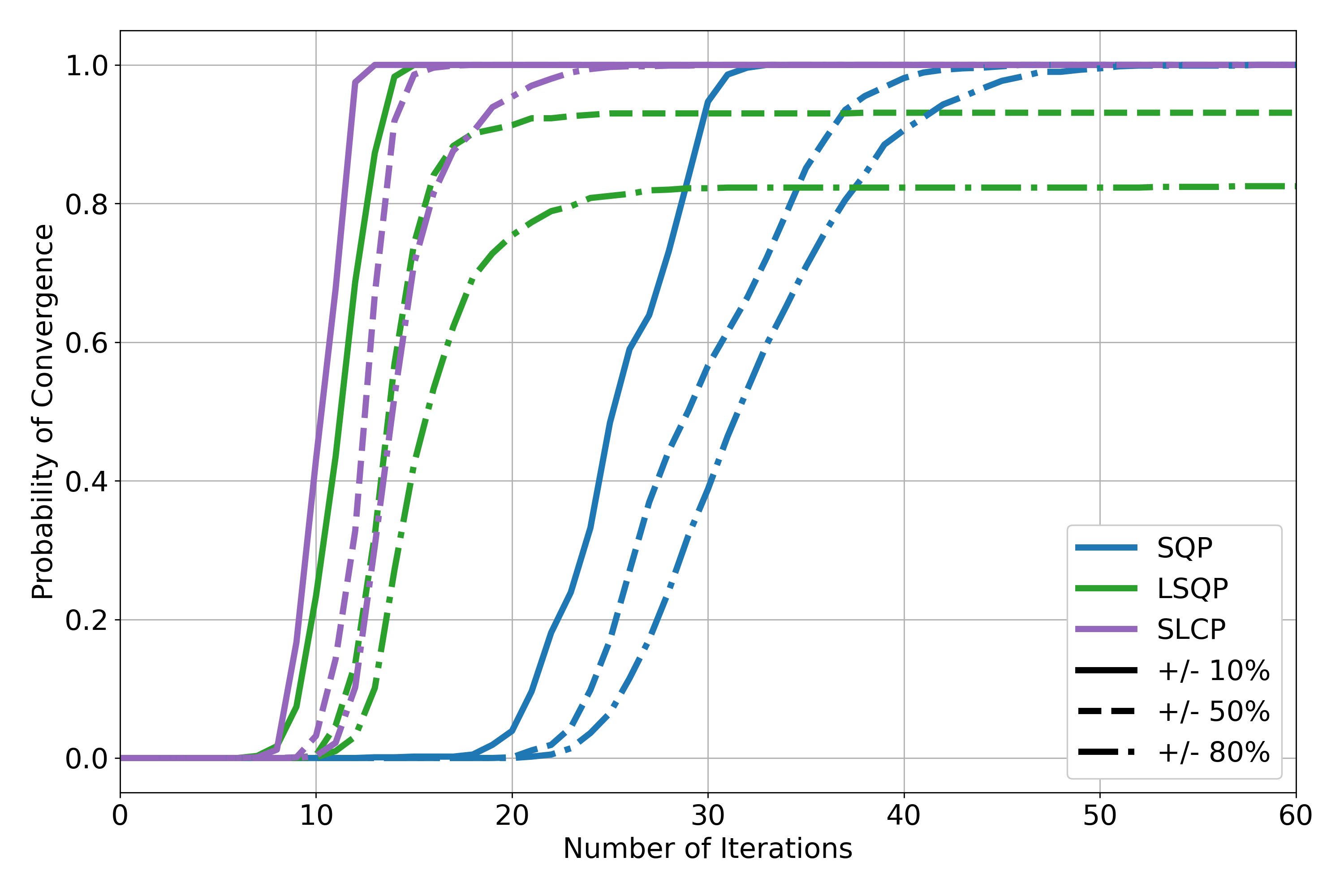} 
\caption{Probability of Convergence vs. Iteration Count for the Floudas Problem}
\label{floudas_results_fig}
\end{figure}

The problem has 8 variables and 6 constraints, none of which are mononomials and only one of which is a posynomial.  The optimal solution is reported by Floudas \cite{floudas2013handbook} for comparison.  Results for this test problem are presented in Figure \ref{floudas_results_fig}, and in Tables \ref{t:benchmarkSummary_good}, \ref{t:benchmarkSummary_reasonable}, and \ref{t:benchmarkSummary_poor}.  

\FloatBarrier
\FloatBarrier

Though this was one of the more interesting case studies for LSQP due to some unexpected tradeoffs with SQP \cite{karcher2021logspace}, the results here paint a clear picture that SLCP strictly outperforms SQP and LSQP for this problem.  
\subsection{Kirschen-Ozturk Problem}
Kirschen \cite{kirschen2018power} proposes the following problem for aircraft sizing using low fidelity analysis models\footnote{Attribution of this problem in the Kirschen paper is given to the uncredited Berk Ozturk \cite{kirschen2018power}, though much of the problem originates in the work of Hoburg \cite{hoburg2014geometric}.}:
\begin{align*}
        \text{minimize}  \quad &  W_f \\
        \text{subject to}  \quad & W_f \geq c_T t D   \\
        & t \geq \frac{R}{V}  \\
        & D \geq \frac12 \rho V^2 S C_D  \\
        & C_D \geq \frac{A_{C_{D_0}}}{S} + k C_f \frac{S_{\text{wet}}}{S} + \frac{C_L^2}{\pi A e}  \\
        & C_f \geq 0.074 Re^{-0.02}  \\
        & Re \leq \frac{\rho V \sqrt{S/A}}{\mu}  \\
        & \frac12 \rho V^2 S C_L \geq W_0 + W_w + \frac12 W_f \\
        & \frac12 \rho V_{\text{min}}^2 S C_{L_{\text{max}}} \geq W \tag{\stepcounter{equation}\theequation}\\\\
        & W \geq W_0 + W_w + W_f  \\
        & W_w \geq W_{w_{\text{surf}}} + W_{w_{\text{strc}}}  \\
        & W_{w_{\text{surf}}} \geq C_{W_{w,1}} S \\
        & W_{w_{\text{strc}}} \geq C_{W_{w,2}} \frac{N_{\text{ult}} A^{\frac32} \sqrt{(W_0 + V_{f_{\text{fuse}}} g \rho_f) W S}}{\tau}  \\
        & V_f \leq V_{f_{\text{avail}}} \\
        & V_f = \frac{W_f}{g \rho_f} \\
        & V_{f_{\text{avail}}} \leq V_{f_{\text{wing}}} + V_{f_{\text{fuse}}} \\
        & V_{f_{\text{wing}}}^2 \leq 0.0009 \frac{S^3}{A} \tau^2\\
        & V_{f_{\text{fuse}}} \leq A_{C_{D_0}} 10[m]  
        \label{kirschenozturk}
\end{align*}
Variables were also constrained to be greater than a small positive constant in order to assist in the construction of sub-problems.

The problem is a signomial program, and is solved by Kirschen using the Difference of Convex Algorithm (DCA) \cite{kirschen2016signomial,burnell2020gpkit,karcher2021logspace}.  That solution is used as the reference solution for this work, but it is important to note that solutions obtained using DCA do not hold any optimality guarantees as the convergence criteria for DCA is based on a relative change in the objective function and not on first order optimality conditions.  

\begin{figure}[htb]
\centering     
\includegraphics[width=0.9\textwidth]{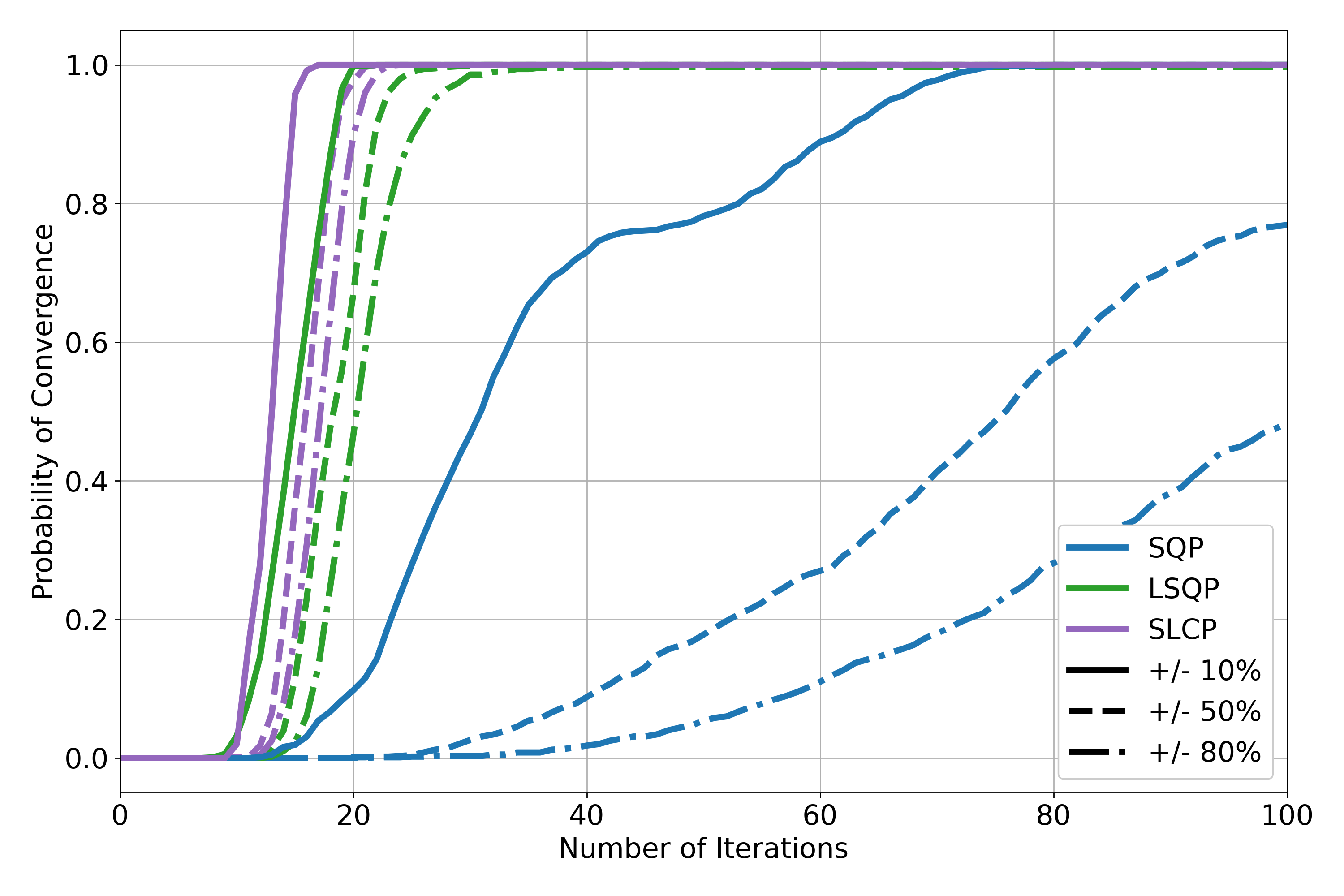} 
\caption{Probability of Convergence vs. Iteration Count for the Kirschen-Ozturk Problem}
\label{ko_results_fig}
\end{figure}

Results for this test problem are presented in Figure \ref{ko_results_fig}, and in Tables \ref{t:benchmarkSummary_good}, \ref{t:benchmarkSummary_reasonable}, and \ref{t:benchmarkSummary_poor}.  

\FloatBarrier
\FloatBarrier

Again, the SLCP algorithm outperforms LSQP, both of which significantly outperform SQP.
\subsection{Defining the Hoburg Problem}
The final test problem comes from Hoburg \cite{hoburg2014geometric}, and is rather extensive, consisting of 82 variables and 119 constraints.  The problem defines the conceptual design and sizing of a UAV, which flies an outbound leg, a return leg, and has a separate set of sprint constraints that are used to size the powerplant.  The problem seeks to minimize the objective:
\begin{equation}
    W_{\text{fuel,out}} + W_{\text{fuel,ret}}
    \label{hoburg_objective}
\end{equation}

\noindent
Subject to the following constraints, which are classified for readability. \\

\noindent
Steady level flight relations:
\begin{equation}
    \begin{aligned}
         W & = \frac12 \rho V^2 C_L S \\
         T & \geq \frac12 \rho V^2 C_D S \\
         Re & = \frac{\rho V S^{1/2}}{A^{1/2} \mu}
    \end{aligned}
    \label{hoburg_slf}
\end{equation}

\noindent
Landing flight condition:
\begin{equation}
    \begin{aligned}
         W_{\text{MTO}} & \leq \frac12 \rho_{\text{sl}} V^2_{\text{stall}} C_{L,\text{max}} S \\
         V_{\text{stall}} & \leq 38 
    \end{aligned}
    \label{hoburg_landing}
\end{equation}

\noindent
Sprint flight condition:
\begin{equation}
    \begin{aligned}
        P_{\text{max}} & \geq \frac{T_{\text{sprint}} V_{\text{sprint}}}{\eta_{0\text{,sprint}}} \\
        V_{\text{sprint}} & \geq 150         
    \end{aligned}
    \label{hoburg_sprint}
\end{equation}

\noindent
Drag model:
\begin{equation}
    \begin{aligned}
         C_D \geq & \frac{0.5}{S} + C_{D_p} + \frac{C_L^2}{\pi e A} \\
         1 \geq & 2.56               \frac{ C_L^{ 5.88}             }{ \tau^{3.32} Re^{1.54} C_{D_p}^{2.26} }+ 
                3.80\times 10^{-9} \frac{ \tau^{ 6.23}            }{ C_L^{0.92} Re^{1.38} C_{D_p}^{9.57}  }+ \\
                & 2.20\times 10^{-3} \frac{ \tau^{ 0.03} Re^{ 0.14} }{ C_L^{0.01} C_{D_p}^{0.73} }+ 
                1.19\times 10^{4}  \frac{C_L^{ 9.78} \tau^{ 1.76} }{ Re^{1.00} C_{D_p}^{0.91} }+ \\
                & 6.14\times 10^{-6} \frac{C_L^{ 6.53}              }{ \tau^{0.52} Re^{0.99} C_{D_p}^{5.19} } 
    \end{aligned}
    \label{hoburg_drag}
\end{equation}

\noindent
Propulsive efficiency:
\begin{equation}
    \begin{aligned}
         \eta_0 & \leq \eta_{\text{eng}} \eta_{\text{prop}} \\
         \eta_{\text{prop}} & \leq \eta_i \eta_v \\
         4 \eta_i + \frac{T \eta_i^2}{\frac12 \rho V^2 A_{\text{prop}} } & \leq 4
    \end{aligned}
    \label{hoburg_prop}
\end{equation}

\noindent
Range constraints:
\begin{equation}
    \begin{aligned}
         R & \geq 5000 \times 10^3 \\
         z_{\text{bre}} & \geq \frac{g R T}{h_{\text{fuel}} \eta_0 W} \\
         \frac{W_{\text{fuel}}}{W} & \geq z_{\text{bre}} + \frac{z_{\text{bre}}^2}{2} + \frac{z_{\text{bre}}^3}{6} + \frac{z_{\text{bre}}^4}{24}
    \end{aligned}
    \label{hoburg_range}
\end{equation}

\noindent
Weight relations:
\begin{equation}
    \begin{aligned}
         W_{\text{pay}} & \geq 500 g \\
         \tilde{W} & \geq W_{\text{fixed}} + W_{\text{pay}} + W_{\text{eng}} \\
         W_{\text{zfw}} & \geq \tilde{W} + W_{\text{wing}} \\
         W_{\text{eng}} & \geq 0.0372 P_{\text{max}}^{0.803} \\
         \frac{W_{\text{wing}}}{f_{\text{wadd}}} & \geq W_{\text{web}} + W_{\text{cap}}\\
         W_{\text{out}} & \geq W_{\text{zfw}} + W_{\text{fuel,ret}} \\
         W_{\text{MTO}} & \geq W_{\text{out}} + W_{\text{fuel,out}} \\
         W_{\text{sprint}} & = W_{\text{out}}
    \end{aligned}
    \label{hoburg_weight}
\end{equation}

\noindent
Wing structural model:
\begin{equation}
    \begin{aligned}
         2q & \geq 1+p \\
         p & \geq 1.9 \\
         \tau & \leq 0.15 \\
         \bar{M}_r & \geq \frac{\tilde{W} A p}{24} \\
         0.92 \bar{w} \tau \bar{t}_{\text{cap}}^2 + \bar{I}_{\text{cap}} & \leq \frac{0.92^2}{2} \bar{w} \tau^2 \bar{t}_{\text{cap}} \\
         8 & \geq \frac{N_{\text{lift}} \bar{M}_r A q^2 \tau}{S \bar{I}_{\text{cap}} \sigma_{\text{max}}} \\
         12 & \geq \frac{A \tilde{W} N_{\text{lift}} q^2}{\tau S \bar{t}_{\text{web}} \sigma_{\text{max,shear}}} \\
         \nu^{3.94} & \geq 0.86 p^{-2.38} + 0.14 p^{0.56}\\
         W_{\text{cap}} & \geq \frac{8 \rho_{\text{cap}} g \tilde{w} \bar{t}_{\text{cap}} S^{3/2} \nu}{3 A^{1/2}} \\
         W_{\text{web}} & \geq \frac{8 \rho_{\text{web}} g r_h \tau \bar{t}_{\text{web}} S^{3/2} \nu}{3 A^{1/2}}
    \end{aligned}
    \label{hoburg_struct}
\end{equation}

\noindent
Additional information is available in the Hoburg paper \cite{hoburg2014geometric}.  This formulation is GP compatible, and therefore has a known global optimum.

Three versions of this problem were considered in an effort to demonstrate the ability of SLCP to systematically evolve model fidelity.  Note that many constraints in the Hoburg formulation are exact, and introduce no uncertainty to the final result (see Equations \ref{hoburg_slf}, \ref{hoburg_landing}, \ref{hoburg_sprint} in particular, though many of the other constraints are also exact).  However, consider the constraint in Equation \ref{hoburg_drag}:

\begin{equation}
         1 \geq 2.56               \frac{ C_L^{ 5.88}             }{ \tau^{3.32} Re^{1.54} C_{D_p}^{2.26} }+ ... +
                6.14\times 10^{-6} \frac{C_L^{ 6.53}              }{ \tau^{0.52} Re^{0.99} C_{D_p}^{5.19} } 
    \label{dragFit}
\end{equation}
\noindent
which is a posynomial fit to a set of XFOIL \cite{drela1989xfoil} data for the NACA 24xx family of airfoils acting as a surrogate model for profile drag coefficient, $C_{D_p}$.  This constraint is enforced for the outbound, return, and sprint segments, and so the form shown in Equation \ref{dragFit} actually represents three separate constraints.  

This model introduces uncertainty in two forms.  First is the uncertainty of the fitted model itself.  For points that were in the original fitting set, the model will not capture the exact $C_{D_p}$ reported by XFOIL because the model fitting process is minimizing some RMS error.  For points not in the original fitting set, interpolation uncertainty is also introduced (ie, features smaller than the sampling interval have been ignored).  This model uncertainty can only be removed by tying XFOIL directly into the optimization problem, which though not possible in Hoburg's original GP formulation can be done with SLCP.  Second is the epistemic uncertainty between the physics modeled in XFOIL and the true underlying physics.  This epistemic uncertainty can only be reduced by using a higher fidelity analysis tool, and so the ability to swap out an integrated XFOIL model with one of higher fidelity would be desirable.  The following case studies will establish a path towards achieving this result.

\subsection{Hoburg Problem as Formulated} \label{hoburgProblem}
\begin{figure}[htb]
\centering     
\includegraphics[width=0.9\textwidth]{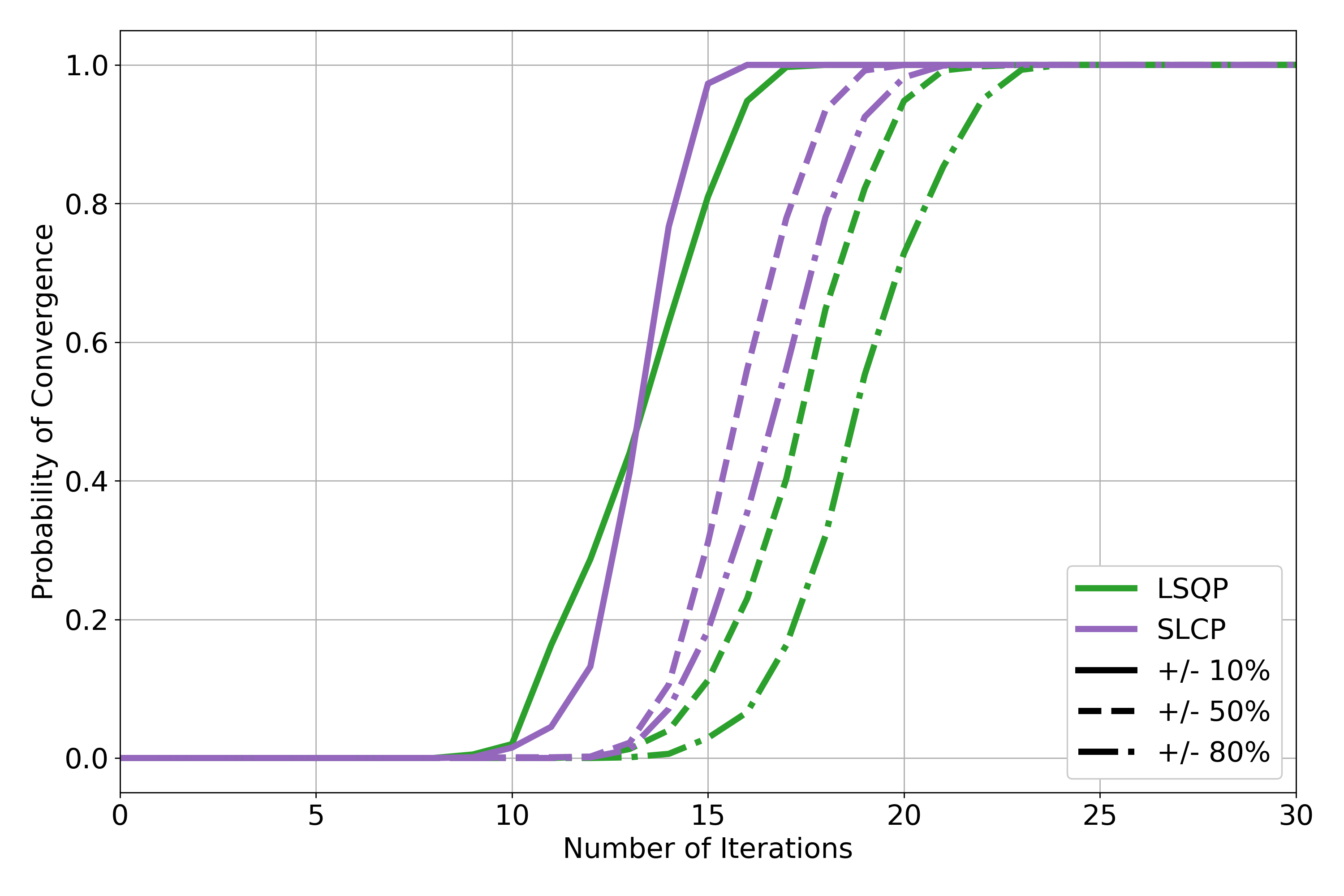}
\caption{Probability of Convergence vs. Iteration Count for the Hoburg Problem with No Black Boxed Constraints}
\label{h0_results_fig}
\end{figure}

First consider solving the Hoburg problem exactly as formulated.  Figure \ref{h0_results_fig} reports the results of the trial runs, with the averages reported in Tables \ref{t:benchmarkSummary_good}, \ref{t:benchmarkSummary_reasonable}, and \ref{t:benchmarkSummary_poor}.  Note that the SQP algorithm was not able to solve this problem as formulated from any initial guess, and so no results are reported for this or any of the subsequent cases.  

\FloatBarrier
\FloatBarrier

There is little difference between the two algorithms when the initial guess is in a region close to the true optimum, since the initialization of the curvature matrix (in this case the identity matrix, see Algorithm \ref{alg:slcp}) dominates the overall performance.  But significant gains are seen as the initial guess decreased in quality, topping out at an average of 11.4\% improvement, or about two iterations.  The significant improvement is due to the posynomial constraints, particularly those characterized by Equation \ref{dragFit}, being exactly represented in the SLCP sub-problem.  The next two cases serve to isolate the amount of computational savings that comes from directly implementing the posynomials of Equation \ref{dragFit}, while also evolving the problem towards something that can be used to integrate higher fidelity models.

\subsection{Hoburg Problem with One Black Boxed Constraint} \label{hoburgProblem1}
Consider that Equation \ref{dragFit} is essentially a representation of the black boxed analysis function:
\begin{equation}
    C_{D_p} = f(C_L, \tau, Re)
    \label{cdpAnalysis}
\end{equation}
imposed in constraint form as:
\begin{equation}
    1 \geq \frac{f(C_L, \tau, Re)}{C_{D_p}}
    \label{cdpCon}
\end{equation}

Many analysis models can be represented by Equation \ref{cdpAnalysis}, including low fidelity methods like XFOIL, high fidelity CFD methods, and even data from wind tunnel tests.  So if the constraint posed in Equation \ref{cdpCon} can be included in the Hoburg formulation, then any analysis model which captures $f(C_L, \tau, Re)$ can be similarly used.

\begin{figure}[htb]
\centering     
\includegraphics[width=0.9\textwidth]{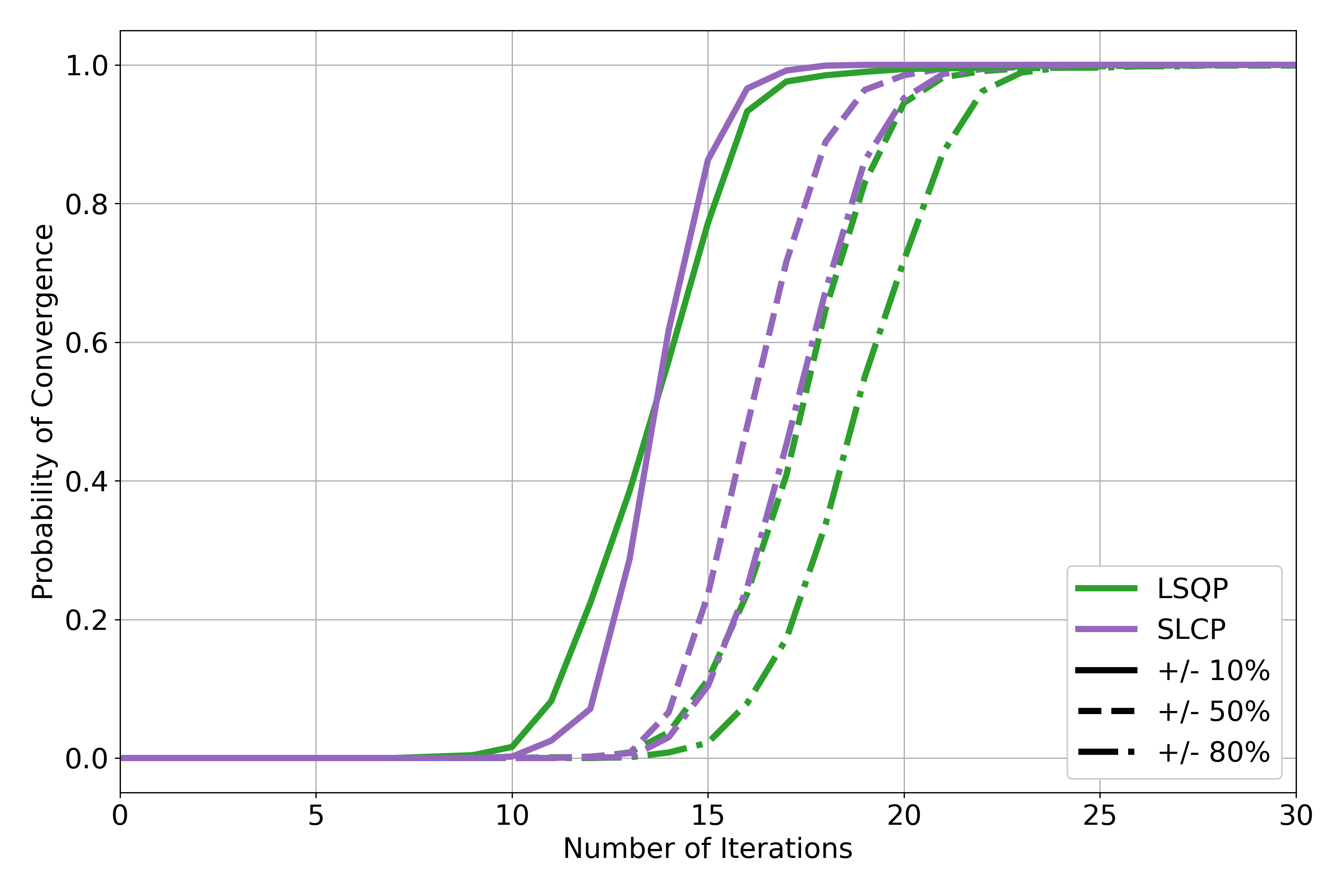} 
\caption{Probability of Convergence vs. Iteration Count for the Hoburg Problem with One Black Boxed Constraint}
\label{h1_results_fig}
\end{figure}

The goal of this test case is to introduce a single black boxed analysis model for $f(C_L, \tau, Re)$ while minimizing the impact on the rest of the GP compatible problem.  Thus a single constraint, in this case the sprint segment instance of Equation \ref{dragFit}, was replaced by Equation \ref{cdpCon}, where $f(C_L, \tau, Re)$ was determined by implicitly solving Equation \ref{dragFit}.  In effect this keeps the function $f(C_L, \tau, Re)$ the same, but hides the true posynomial from the SLCP algorithm and should therefore reduce computational efficiency when compared to the previous case while maintaining the same solution.  Results are reported in Figure \ref{h1_results_fig} and Tables \ref{t:benchmarkSummary_good}, \ref{t:benchmarkSummary_reasonable}, and \ref{t:benchmarkSummary_poor}.  

\FloatBarrier
\FloatBarrier

As expected, the loss of a posynomial constraint reduces the efficiency of SLCP and moves these curves in Figure \ref{h1_results_fig} closer to the LSQP curves, however SLCP still demonstrates a performance gain between 6-8\% as the initial guesses get worse.

\subsection{Hoburg Problem with Three Black Boxed Constraints} \label{hoburgProblem3}
\begin{figure}[htb]
\centering     
\includegraphics[width=0.9\textwidth]{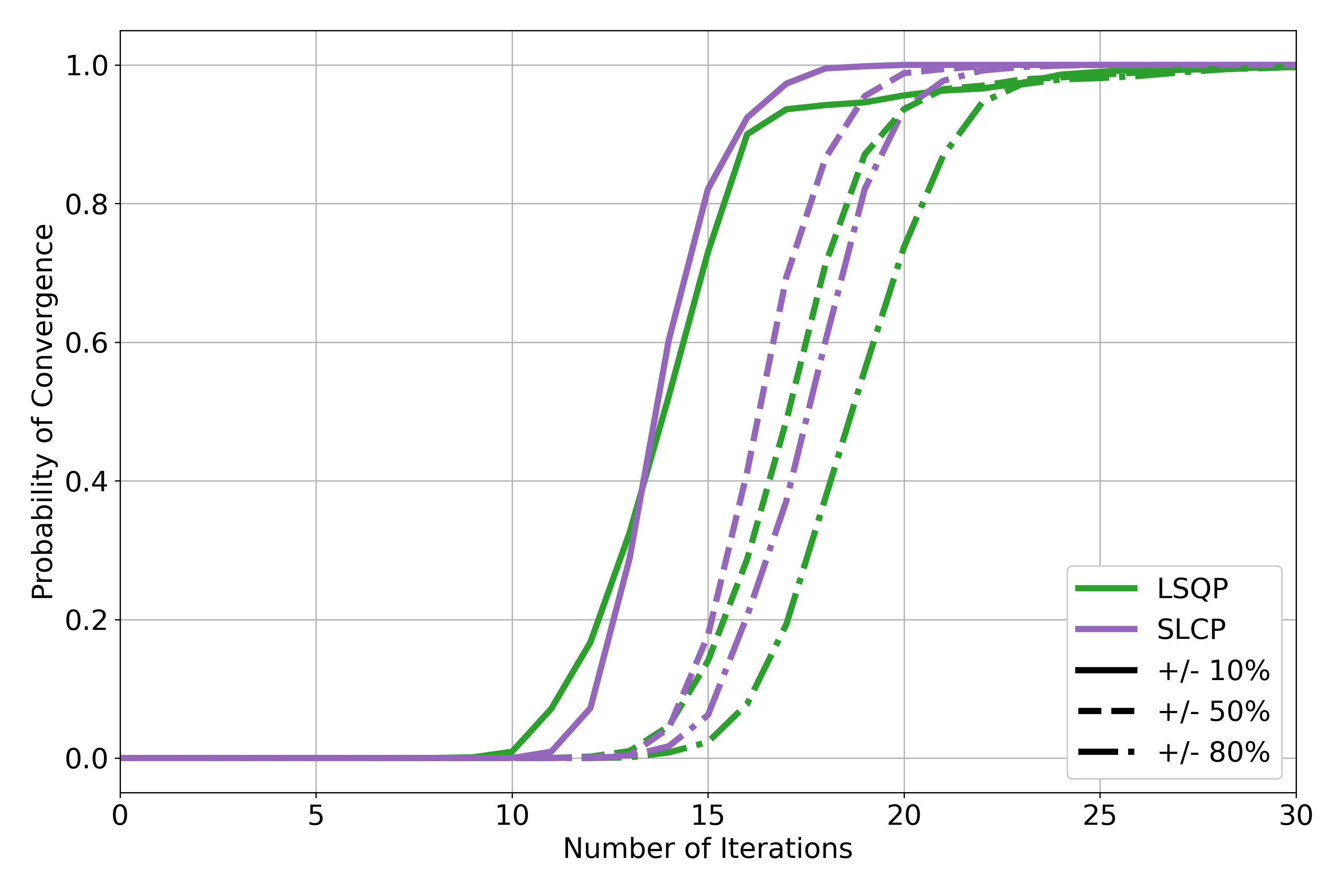} 
\caption{Probability of Convergence vs. Iteration Count for the Hoburg Problem with Three Black Boxed Constraints}
\label{h3_results_fig}
\end{figure}

For the final test problem, all three constraints represented by Equation \ref{dragFit} (outbound, return, and sprint) were replaced with black boxed models for $f(C_L, \tau, Re)$.  As with the previous case, the black boxes were implicit implementations of Equation \ref{dragFit}, leaving the problem identical to the original Hoburg problem but with these three posynomials hidden to the SLCP algorithm.  The results are reported in Figure \ref{h3_results_fig} and in Tables \ref{t:benchmarkSummary_good}, \ref{t:benchmarkSummary_reasonable}, and \ref{t:benchmarkSummary_poor}

\FloatBarrier
\FloatBarrier

As expected, the elimination of two more posynomial constraints further shrinks the gap between the SLCP and LSQP curves in Figure \ref{h3_results_fig}.  But critically, the fact that all three constraints have been black boxed means that the posynomial surrogate from Equation \ref{dragFit} has been entirely eliminated, and a higher fidelity analysis model can simply be swapped in as a new black box that represents $f(C_L, \tau, Re)$.

To demonstrate this procedure, XFOIL was directly tied in to the optimization formulation, with gradients being computed using a finite difference scheme.  Note that unlike the previous three problem formulations, the inclusion of XFOIL directly has fundamentally altered the structure of the problem and will affect the result (see Table \ref{t:hoburg_x} for a summary of the changes to the optimal solution).  The GP optimal solution was used as the initial guess, and SLCP reached the optimal soluiton in 32 iterations.    

\begin{table}[htbp] 
  \footnotesize 
  \begin{center} 
  \caption{Comparison of variables with a greater than 1\% variation when XFOIL is tied directly into the optimization formulation} 
  \label{t:hoburg_x} 
  \begin{tabular}{ c  c  c  c  } 
      \hline
        \rule{0pt}{2.25ex} Design Variable    & GP Optimal Value & SLCP+XFOIL Optimal Value & Percent Change  \\ \hline \rule{0pt}{2.25ex} 
        $C_{D_p\text{,sprint}}$         &  0.005732         &   0.005455        &  -4.82  \% \\
        $C_{D\text{,sprint}}$           &  0.007760         &   0.007477        &  -3.65  \% \\
        $P_{\text{max}}$ [kW]           &  1186.1           &   1143.2          &  -3.61  \% \\
        $T_{\text{sprint}}$ [N]         &  2209.6           &   2133.9          &  -3.42  \% \\
        $C_{D_p\text{,out}}$            &  0.005455         &   0.005615        &   2.93  \% \\
        $W_{\text{eng}}$ [N]            &  2805.8           &   2724.1          &  -2.91  \% \\
        $C_{D_p\text{,ret}}$            &  0.005470         &   0.005627        &   2.87  \% \\
        $C_{D\text{,ret}}$              &  0.012682         &   0.012894        &   1.67  \% \\
        $C_{D\text{,out}}$              &  0.012668         &   0.012876        &   1.64  \% \\
        $\bar{t}_{\text{cap}}$          &  0.004544         &   0.004611        &   1.47  \% \\
        $\bar{I}_{\text{cap}}$          &  2.021e-5         &   2.555e-5        &   1.37  \% \\
        $W_{\text{cap}}$ [N]            &  4103.8           &   4158.4          &   1.33  \% \\
        $W_{\text{wing}}$ [N]           &  8613.4           &   8723.7          &   1.28  \% \\
        $C_{Di\text{,sprint}}$          &  0.000232         &   0.000230        &  -1.12  \% \\
        $C_{Di\text{,ret}}$             &  0.005416         &   0.005475        &   1.08  \% \\
        $C_{L\text{,ret}}$              &  0.574            &   0.580           &   1.03  \% \\
        $W_{\text{fuel, out}}$ [N]      &  3082.6           &   3114.1          &   1.02  \% \\  \hline
  \end{tabular} 
 \end{center} 
\end{table} 
\FloatBarrier

Since XFOIL has been entirely contained in a black box of $f(C_L, \tau, Re)$, higher fidelity tools like MSES \cite{drela2007user}, 2D Euler methods, or 2D Navier-Stokes simulations can readily integrated in exactly the same fashion.  Furthermore, the design variables available in these higher fidelity analysis models can be optimized as a part of the SLCP run.  For example, a NACA 4-series airfoil could be parameterized as $f(C_L, \tau, Re, c_{max}, l_{c_{max}})$ where $c_{max}$ is the maximum camber in fraction of chord and $l_{c_{max}}$ is the location of that maximum camber, again as a fraction of chord.  In this way, geometry fidelity can also evolve with the evolutions in analysis fidelity, all while anchored in a strong initial guess migrated from the previous fidelity level.

Together, Geometric Programming and Sequential Log-Convex Programming offer a potential approach to engineering design when multiple levels of analysis fidelity must be considered:
\begin{enumerate}
    \item Develop low-fidelity conceptual models in a GP compatible form
    \item Solve the GP to obtain a low fidelity candidate design to be used as an initial guess for higher levels of fidelity
    \item Evolve analysis models to higher fidelity as desired and implement them in the optimization formulation as black box functions
    \item Add any new design variables made available by the higer fidelity analysis models to the optimization formulation
    \item Solve the new optimization problem with SLCP, using the GP solution as an initial guess
\end{enumerate}

Future work will dive into this methodology in great depth, particularly in developing a systematic process for evolving analysis models and the integration of new design variables into the outer optimization loop.

\begin{table}[htbp] 
  \footnotesize 
  \begin{center} 
  \caption{Average over 1000 trials of the number of iterations required to obtain an optimal solution starting from a guess within +/- 10\% of the known optimum} 
  \label{t:benchmarkSummary_good} 
  \begin{tabular}{c  c  c  c } 
      \cline{2-4}
      \rule{0pt}{2.25ex} & SQP   & LSQP             & SLCP             \\ \hline \rule{0pt}{2.25ex}
      Floudas            & 25.87 & 11.70 (-54.77\%) & 10.74 (-58.48\%) \\ 
      Kirschen-Ozturk    & 33.93 & 15.36 (-54.73\%) & 13.34 (-60.68\%) \\ 
      Hoburg-0           & -     & 13.70            & 13.65 (-0.32\%)  \\ 
      Hoburg-1           & -     & 14.06            & 14.17 (+0.78\%)  \\ 
      Hoburg-3           & -     & 14.65            & 14.32 (-2.27\%)  \\ \hline 
  \end{tabular} 
 \end{center} 
\end{table} 

\begin{table}[htbp] 
  \footnotesize 
  \begin{center} 
  \caption{Average over 1000 trials of the number of iterations required to obtain an optimal solution starting from a guess within +/- 50\% of the known optimum} 
  \label{t:benchmarkSummary_reasonable} 
  \begin{tabular}{ c  c  c  c  } 
      \cline{2-4}
      \rule{0pt}{2.25ex} & SQP   & LSQP             & SLCP             \\ \hline \rule{0pt}{2.25ex}
      Floudas            & 30.17 & 14.30 (-52.60\%) & 12.92 (-57.18\%) \\ 
      Kirschen-Ozturk    & 72.55 & 18.89 (-73.96\%) & 16.38 (-77.42\%) \\ 
      Hoburg-0           & -     & 17.79            & 16.29 (-8.44\%)  \\ 
      Hoburg-1           & -     & 17.82            & 16.66 (-6.47\%)  \\ 
      Hoburg-3           & -     & 17.66            & 16.87 (-4.49\%)  \\ \hline 
  \end{tabular} 
 \end{center} 
\end{table} 

\begin{table}[htbp] 
  \footnotesize 
  \begin{center} 
  \caption{Average over 1000 trials of the number of iterations required to obtain an optimal solution starting from a guess within +/- 80\% of the known optimum} 
  \label{t:benchmarkSummary_poor} 
  \begin{tabular}{ c  c  c  c  } 
      \cline{2-4}
      \rule{0pt}{2.25ex} & SQP   & LSQP             & SLCP              \\ \hline \rule{0pt}{2.25ex}
      Floudas            & 32.79 & 16.29 (-50.32\%) & 14.92 (-54.50\%)  \\ 
      Kirschen-Ozturk    & 91.18 & 21.06 (-76.90\%) & 17.66 (-80.63\%)  \\ 
      Hoburg-0           & -     & 19.33            & 17.13 (-11.40\%)  \\ 
      Hoburg-1           & -     & 19.30            & 17.70 (-8.30\%)   \\ 
      Hoburg-3           & -     & 19.33            & 18.02 (-6.77\%)   \\ \hline 
  \end{tabular} 
 \end{center} 
\end{table} 
\FloatBarrier

\section{Conclusions}
The method of Sequential Log-Convex Programming proposed in this work brings together the efficiency of Geometric Programming and the flexibility of Sequential Quadratic Programming in a practical algorithm for engineering design.  An engineer who begins the design process using simple GP compatible models can utilize SLCP to evolve analysis fidelity and reduce uncertainty, while an engineer who begins the design process with higher fidelity models can utilize SLCP to exploit underlying mathematical structure inherent in many design problems and reduce overall computational expense with little to no modification of existing design practices.  

\begin{acknowledgements}
This material is based on research sponsored by the U.S. Air Force under agreement number FA8650-20-2-2002. The U.S. Government is authorized to reproduce and distribute reprints for Governmental purposes notwithstanding any copyright notation thereon.  The views and conclusions contained herein are those of the authors and should not be interpreted as necessarily representing the official policies or endorsements, either expressed or implied, of the U.S. Air Force or the U.S. Government.

The authors would like to thank Mark Drela, Marshall Galbraith, John Hansman, and John Dannenhoffer for their input into the technical matter, along with the EnCAPS Technical Monitor Ryan Durscher.

The authors also acknowledge the MIT SuperCloud and Lincoln Laboratory Supercomputing Center for providing high performance computing, database, and consultation resources that have contributed to the research results reported within this paper.

\end{acknowledgements}
\bibliographystyle{spmpsci}      
\bibliography{bib_deck}   

\end{document}